\documentclass[12pt]{article}

\usepackage{amsmath,amsthm,amssymb}

\theoremstyle{plain}
\newtheorem{Thm}{Theorem}
\newtheorem{Prop}{Proposition}[section]
\newtheorem{Cor}{Corollary}[section]
\newtheorem{lem}{Lemma}[section]

\makeatletter

\@addtoreset{equation}{section}

\def\R{\mathbb{R}}
\def\Z{\mathbb{Z}}

\def\v2{\vskip2mm}
\def\n{\noindent}

\def\({(\!(}
\def\){)\!)}
\def\F{{\cal F}}
\def\a{\alpha}

\def\e{\varepsilon}
\def\de{\delta}

\def\la{\lambda}

\def\th{\theta}

\def\Ga{\Gamma}
\def\La{\Lambda}
\def\Om{\Omega}

\def\pf{{\it Proof.}}
\def\v2{\vskip2mm}
\def\n{\noindent}

\def\0{{\bf 0}}

\def\tst12{{\textstyle \frac12}}
\def\1{{\bf 1}}

\def\n{\noindent}
\def\beq{\begin{eqnarray*}}
\def\eeq{\end{eqnarray*}}

\def\beqn{\begin{equation}}
\def\eeqn{\end{equation}}

\begin{document}

\begin{center}
{\bf Scaling  limits of    random walk bridges conditioned to  avoid a finite set} \\
\vskip4mm
{K\^ohei UCHIYAMA} \\
\vskip2mm
{Department of Mathematics, Tokyo Institute of Technology} \\
{Oh-okayama, Meguro Tokyo 152-8551\\
}
\end{center}

\vskip8mm
\n
\vskip8mm
\n
{\it running head}:    Scaling  limits of random walk bridges 

\vskip2mm
\n
{\it key words}:  
 random walk on the integer lattice;  conditioned to  avoid a  set; third moment;  functional limit theorem; killing on a finite set; tightness of pinned walk; tunneling. 
\vskip2mm
\n

\vskip6mm 

\begin{abstract}
This paper concerns  a scaling limit of  a  one-dimensional random walk $S^x_n$ started from $x$  on the integer lattice conditioned to avoid a non-empty finite set $A$,  the random walk  being assumed to be  irreducible and have zero mean. Suppose the  variance $\sigma^2$ of the increment law  is finite. Given positive constants $b$, $c$ and $T$ we consider the scaled process $S^{b_N}_{[tN]}/\sigma\sqrt N$, $0\leq t \leq T$  started from a point $b_N \approx b\sqrt N$ conditioned to arrive at another point 
$\approx -c\sqrt N$ at $t=T$ and avoid $A$ in between and discuss  the functional  limit of it as $N\to\infty$.  We show that it converges in law to a continuous process 
if  $E[|S_1|^3; S_1<0] <\infty$. If $E[|S_1|^3; S_1<0] =\infty$ we suppose   $P[S_1<u]$ to vary regularly as $u\to -\infty$ with exponent $-\beta$, $2\leq \beta\leq 3$ and show that  it converges to a process which  has one downward  jump  that clears the origin  if $\beta<3$; in  case  $\beta=3$ there arises the same  limit process as  in  case  $E[|S_1|^3; S_1<0] <\infty$. In  case $\sigma^2=\infty$  we  consider the  special case  when  $S_1$ belongs  to the domain of   attraction of a stable law of index $1<\alpha <2$ having no negative jumps  and obtain   analogous results. %

\end{abstract}
\vskip6mm

\vskip6mm
\section{Introduction}
Let $S_n, n=0, 1, 2\ldots$ be a random walk on $\Z$ started at the origin, namely $S_0=0$ and
 $S_{n+1}- S_{n}$ are i.i.d. random variables taking values in the integer lattice  $\Z$. Let   $(S_n)$ be defined on a probability space $(\Om, \F, P)$ and suppose that  $(S_n)$ is irreducible and 
 $$E S_1= 0.$$
Let $\sigma$ be the variance of the step variable:
$\sigma =(E[|S_1|^2])^{1/2}$.  We consider the both of cases $\sigma <\infty$ and  $\sigma=\infty$, but 
usually suppose $\sigma<\infty$  unless the contrary is stated  explicitly when we discuss the problem for  the case $\sigma=\infty$.
In  \cite{U1dm.f.s} the present author obtained a precise asymptotic form of transition probability of  the walk $S$ killed on a finite non-empty set $A$ (in case  $\sigma<\infty$). In the present paper we are interested in  the behaviour  of $S^{b_N}_k:= S_k + b_N, k=1,\ldots, N$, for large $N$, conditioned on the events  $ S^{b_N}_k \notin A$  for $k\leq N$ and $ S^{b_N}_N= - b_N$ where  $b_N \approx b\sqrt N$, $b>0$. In \cite{U1dm.f.s} it is observed that  if $E[|S_1|^3; S_1<0]<\infty$, then the walk thus conditioned  \lq\lq continuously'' transits from the positive half line to the negative half still avoiding $A$, while  if $E[|S_1|^3; S_1<0] =\infty$  it clears $A$ by one \lq\lq long jump''. In this paper we observe that  this is reflected to  the scaling limit.  In case $E[|S_1|^3; S_1<0] =\infty$ we suppose that  $P[S_1<u]$ is regularly varying as $u\to -\infty$ with index $-\beta$, $2\leq \beta\leq3$.  Then we prove  that  the scaled process converges in law  and the limit process is continuous in the former case; in the latter case   it has exactly one downward  jump if $\beta<3$ while  the limit process agrees with that  of the former case if $\beta=3$. In case $\sigma =\infty$    analogous results are given for the  special case when  $S_1$ belongs  to the domain of  attraction of a stable law with exponent $1<\a<2$  having no negative jumps.

There are  a lot of  works dealing with various problems concerning   random walks on the real line conditioned to avoid a finite set or a half line. To mention among them only those studying the functional limit theorems,  
  the finite set case  are studied by  \cite{B2} and \cite{Kai}, whereas  for the half line we have a long list of  papers for which  the readers are referred to   \cite{Bt-D}, \cite{Car} where brief descriptions of them are  found. 

\section{Statements of results}

Supposing $\sigma<\infty$  we  first describe the processes arising in the limit of the random  walk bridges  conditioned to avoid  $A$,
and then state the convergence results. The case $\sigma=\infty$ will be discussed separately after that. There appear some processes in the limit, which we suppose to be given in the same probability space as $(S_n)$.
\v2
 {\bf Limit processes.} Let $W_t,t\geq 0$ be a standard linear Brownian motion started at $0$,
   $Y^x_t$, $x, t\geq 0$ a 3-dimensional Bessel process started at $x$, and  $Y^{x, y,t}_s$ ($ 0\leq s\leq t$,  $x\geq 0, y\geq 0$) the Bessel bridge of length $t$ joining 
 $x$ and  $y $ obtained  by conditioning $Y^x$.  The processes   $Y^{x}, x \geq 0$ are supposed to be  independent  
of $W$. For $x\in \R$ we write  $W^x_t$ for $W_t+x$ and let $\sigma_0^{x, W}$ be the first passage
time of zero by $(W^x_t)$.
 
 Put
$$\mathfrak{g}_t(x) =\frac1{\sqrt{2\pi t}} e^{-x^2/2t}, \quad \mathfrak{g}^0_t(x,y) =\mathfrak{g}_t(y-x) - \mathfrak{g}_t(x+y)$$
and
$$\quad  \rho_t(x) = \frac{|x|}{t}\mathfrak{g}_t(x).$$
For each constants $T>0$ and $c>0$,  $Y^{0,c,T}_t, 0\leq t\leq T$ is a Markov process with transition law
 \beqn\label{trs_d}
\begin{array}{cl}
{\displaystyle \frac{P[Y^{0,c,T}_t\in dz]}{dz} = \frac{\rho_t(z)\mathfrak{g}^0_{T-t}(z,c)}{
\rho_T(c)},} \\[3mm]
{\displaystyle \frac{P[Y^{0, c,T}_t\in dz\,|\, Y^{0, c,T}_s=w]}{dz} =\frac{ \mathfrak{g}^0_{t-s}(w,z) 
\mathfrak{g}^0_{T-t}(z,c)}{ \mathfrak{g}^0_{T-s}(w,c)}}
\end{array}
\eeqn
(cf. \cite[Section VI.3]{RY}). In the sequel the letters $b, c$ and $T$  denote positive constants.
\v2 
 
{\sc  Transit made by creeping. } \, Writing $\tau = \sigma^{b,W}_0 $ we  define the time-inhomogeneous process $X^{b,c, T}_t$, $0\leq t\leq T$ by
 \beqn\label{X0}
 X^{b,c, T}_t =\left\{ \begin{array} {ll} W^b_t \quad & 0\leq t\leq \sigma^{b,W}_0,\\
 - Y^{0, c,\, T-\tau} _{t-\tau} \quad & \tau = \sigma_0^{b,W} < t\leq T,
 \end{array} \right.
 \eeqn
on the event   $\tau <T$.  The process $X^{b,c, T}$ conditioned on 
 $\{ \tau <T\}$ is a  Markov process on $\R$ (inhomogeneous in time) whose transition probability density  $q(s,x; t,y)$  is described  below.  (See (\ref{T00}) for the finite dimensional distribution.)   Note that   if $\tilde W$ is a linear Brownian motion independent of $W$, the Bessel bridge $Y^{0, c,\, T-\tau} _{t-\tau} $ in (\ref{X0}) can be substituted for by   
 $\tilde W^c_{T- t}$ conditional on $\sigma^{c,\tilde W}_0= T-\tau$, the two processes having the same law (cf. e.g., \cite[XI(3.12)]{RY}), provided the event $\{ \tau <T\}$ conditioning  $X^{b,c, T}_t $ is accordingly replaced by $\tau + \sigma^{c,\tilde W}_0 =T$.

 For $0<s <t<T$,
\beqn\label{m1}
q(s,x; t,y) =\left\{ \begin{array} {ll} {\displaystyle  \frac{\mathfrak{g}^0_{t-s}(x,y)\rho_{T-t}(y+c)}{\rho_{T-s}(x+c)} \quad} &   x\geq 0,\, y > 0 ,\\[4mm]{\displaystyle 
\frac{\rho_{t-s}(x-y)\mathfrak{g}^0_{T-t}(y,-c)}{\rho_{T-s}(x+c)}} \quad& y\leq 0\leq x,\\[4mm]{\displaystyle 
\frac{\mathfrak{g}^0_{t-s}(x,y)\mathfrak{g}^0_{T-t}(y,-c)}{\mathfrak{g}^0_{T-s}(x,-c)}} &   x < 0,\, y\leq 0, \\[4mm]
0       &  x < 0 < y,
\end{array} \right.
 \eeqn
and for $0< t<T$,
\beqn\label{m2}
q(0,b; t,x) =q(t, x; T,-c) = \left\{ \begin{array} {ll} {\displaystyle  \frac{\mathfrak{g}^0_{t}(b,x)\rho_{T-t}(x+c)}{\rho_{T}(b+c)} \quad} & 0\leq x,\\[4mm]{\displaystyle 
\frac{\rho_{t}(b-x)\mathfrak{g}^0_{T-t}(x,-c)}{\rho_{T}(b+c)}} \quad& x <0.\\
\end{array} \right.
 \eeqn
\v2
{\sc Remark 1.} \, (a)\, The semi group property of 
$q(s,x; t,y)$ can be directly ascertained 
by using the relations
$$\int_0^\infty \rho_t(z)\mathfrak{g}^0_s(z,y)dz = \rho_{t+s}(y) \quad \mbox{and}\quad \int_0^t \rho_{t-s}(x)\rho_s(y)ds = \rho_t(x+y) \qquad (x, y>0).$$
(e.g.,  to see $\int_{-\infty}^\infty q(0,b;t,x)dx=1$  write $\rho_t(b-x)=\int_0^t\rho_{t-s}(-x) \rho_s(b)ds$ and similarly for $\rho_{T-t}(x+c)$).
The former one says that the family $(\rho_t(z)dz)_{ t>0, x>0}$ constitutes  the entrance law for the semigroup $\mathfrak{g}^0_t(x,y)dy$ (cf., e.g., \cite{RY}) and the latter follows from the fact that $\rho_\cdot(y-x)$, $y>x$ is the transition density of the passage time process $(\sigma^{0,W}_x)_{x\geq 0}$ (\cite{IM}). 

(b) \, The process  $X^{b,c, T}$ defined above can be obtained as a normalized  limit of Brownian bridge killed at   rate $\la \ell(t)$ as $\la\to\infty$, where $\ell(t)$ is the Brownian  local time at zero  (see Appendix (A)).

(c) \, Although the limit processes are described by means of the 3-dimensional Bessel process, one may think  that there should naturally appear the Brownian meander.    The bridges of  the two processes have  the same law which  can be described by the Brownian motion killed on hitting  zero together with the entrance law for it  as is displayed in  (\ref{trs_d}).   (See also Remark 3.)

\vskip4mm
{\sc Transit  made by a jump. } \, Let $2\leq \beta <3$. Put for $t>0$ and $ y<0<x$,
$$J_{t}(x,y) = \int_0^t ds \int_0^\infty \mathfrak{g}_s^0(x,w)dw\int_0^\infty \frac{\beta}{(w+z)^{1+\beta}} \mathfrak{g}_{t-s}^0(-z,y)dz$$
(see Remark 2 below for the integrability).  
 Let $\tau$ and $\eta$  be positive random variables dependent on $W^b$  with the conditional law
\[
\begin{array}{rr}
{\displaystyle P[\tau\in dt, \eta \in dz\,|\, \F_t^{W,b}] = \frac{\beta\mathfrak{g}_{T-t}^0(-z,-c)}{(W^b_t+z)^{1+\beta}}\cdot \frac{\1(\sigma_0^{b,W}>t) dt dz}{J_T(b,-c)},}&\\[4mm]
  0<t< T, z>0.& 
\end{array}
\]
where $\F_t^{W,b}$ is the $\sigma$-field on $\Om$ generated by  $W^b_s, s\leq t$. Note that 
$$E\bigg[\int_0^Tdt\int_0^\infty \frac{\beta {\mathfrak{g}}_{T-t}^0(-z,-c)}{(W^b_t +z)^{1+\beta}}\1(\sigma_0^{b,W}>t) dz \bigg]\frac{1}{J_T(b,-c)} =1$$
and $\tau$ is an $\F^{W,b}_t$-stopping time taking values in $(0,T)$ a.s.
 Let the Bessel processes   $Y^{x}$, $x \geq 0$ be
  independent of $(\tau, \eta)$ as well as of  $W$.  Define
\beqn\label{X1}
\breve X^{b,c, T}_t =\left\{ \begin{array} {ll} W^b_t \quad & 0\leq t < \tau,\\
- Y^{\eta, c,\,T-\tau}_{t-\tau} \quad & \tau \leq t< T.
 \end{array} \right.
 \eeqn
 Then  $\breve X^{b,c, T}_t$, $t\leq T$  is a  Markov process on $\R\setminus \{0\}$ and  its transition probability is given by replacing  $\rho$ by $J$ in (\ref{m1}) and (\ref{m2}):
For $0<s <t<T$,
$$\breve q(s,x; t,y) =\left\{ \begin{array} {ll} 
{\displaystyle  \frac{\mathfrak{g}^0_{t-s}(x,y)J_{T-t}(y,-c)}{J_{T-s}(x,-c)} \quad} &    x\geq 0, \, y>0 \\[4mm]{\displaystyle 
\frac{J_{t-s}(x,y)\mathfrak{g}^0_{T-t}(y,-c)}{J_{T-s}(x,-c)}} \quad& y\leq 0\leq x,\\[4mm]{\displaystyle 
\frac{\mathfrak{g}^0_{t-s}(x,y)\mathfrak{g}^0_{T-t}(y,-c)}{\mathfrak{g}^0_{T-s}(x,-c)}} &   x < 0,\, y\leq 0, \\[4mm]
0       &  x < 0 < y,
\end{array} \right.
 $$
and for $0< t<T$,
$$\breve q(0,b; t,x) = \breve q(t, x; T,-c) = \left\{ \begin{array} {ll} {\displaystyle  \frac{\mathfrak{g}^0_{t}(b,x)J_{T-t}(x,-c)}{J_{T}(b,-c)} \quad} & 0\leq x,\\[4mm]{\displaystyle 
\frac{J_{t}(b,x)\mathfrak{g}^0_{T-t}(x,-c)}{J_{T}(b,-c)}} \quad& x <0.\\
\end{array} \right.
 $$
 
 \v2
{\sc Remark 2.} The integrability of the repeated  integral defining $J_t$ reduces to that of  $\int_0^1dz\int_0^1zw(z+w)^{-1-\beta}dw$ which is finite if and only if $\beta<3$.

\vskip4mm
{\bf Scaling limits.} Let $(S_n)$ be the random walk on  $\Z$ as specified in Introduction.   Let
$$S^x_n = x+S_n\quad (x\in \Z)$$
and write $\sigma^x_D$ for the first hitting time of a non-empty set $D\subset \Z$ by $S^x$:
$$\sigma^x_D =\inf \{n\geq 1: S^x_n\in D\}.$$
We shall write   $S^x_{\sigma_{\! D}}$ for  $S^x_{\sigma^x_D}$. 
Let  $b$ and $c$ be positive constants,  and $b_N$ and $c_N$, $N=1,2,\ldots$ two sequences of integers   such that as $N\to\infty$
\beqn\label{b/c}
b_N \sim b \sqrt{ \sigma^2N} 
\quad \mbox{and} \quad c_N\sim c\sqrt{ \sigma^2 N} 
\eeqn
(provided   $\sigma^2 <\infty$),  and define the scaled process
\beqn\label{X} X^{(N)}_\cdot = \mbox{the linear interpolation of}\; (S^{b_N} _{k}/\sqrt{\sigma^2 N}: k=0,1, 2,\ldots  ),
\eeqn  
where $\sim$ means that the ratio of two sides of it  approaches  1. Let  $T$ be an arbitrarily given positive constant and $A$ a non-empty finite subset of $\Z$.  We are concerned with   the law of $X^{(N)}$ under the condition that  
\beqn\label{CD}
 \sigma_A^{b_N}>NT \quad\mbox{ and}\quad   S^{b_N}_{\lfloor NT\rfloor}= -c_N,
 \eeqn
namely  the walk $S^{b_N} $ visits $-c_N$ at time $\lfloor NT\rfloor$ without entering  $A$
($\lfloor \cdot \rfloor$ stands for the integer part). 
We suppose that $P[ \sigma^x_{\{y\}} <\infty]>0$ for some, then all, $x>\max A$ and $y<\min A$ so that the probability of event (\ref{CD}) is positive for all sufficiently large $N$  if the walk is further supposed to be  temporally aperiodic (which  is not restrictive for the present problem).  In what follows  the probability of (2.8) is tacitly  supposed to be positive  when the conditioning on (2.8) is considered. 

Let $F$ be the distribution function of $S_1$: $F(u)=P[S_1\leq u]$, $u\in \R$. 

\v2
\begin{Thm}\label{thm1} \, Suppose  that  $E[|S_1|^3; S_1<0]<\infty$ or  $F(u)$ is regularly varying with index $-3$ as $u\to -\infty$.  Then the  law of $( X^{(N)}_t: 0\leq t\leq T) $ conditioned  on   event (\ref{CD})
  converges to the law of $(X^{b,-c,T}_t)_{0\leq t\leq T}$ conditioned on $\sigma^{b,W}_0<T$ relative  to the uniform topology of $C([0,T],\R)$.
\end{Thm}

\begin{Thm} \label{thm2} \, Suppose that    $F(u)$ is regularly varying with index $-\beta$ as $u\to -\infty$ ,  $2\leq \beta<3$.  
     Then the   law of $(X^{(N)}_t: 0\leq t\leq T) $ conditioned  on   event (\ref{CD})  converges to the law of $(\breve X^{b,-c,T}_t)_{0\leq t\leq T}$
relative to the Skorohod topology of $D([0,T], \R)$.
\end{Thm}
\v2

Here we briefly explain how  the two types of limit processes emerge in Theorems \ref{thm1} and \ref{thm2} and thereby indicate  a crucial point in question. 
For $N$ large  the way   $S^{b_N}$  enters the half line $(-\infty,0]$ is virtually determined by
 \beqn\label{H+}
 H_{(-\infty, 0]}^{+\infty}(y) = \lim_{x\to\infty} P[S^x_{\sigma_{(-\infty,0]}} =y],\quad y\leq 0,
 \eeqn 
 the hitting distribution of $(-\infty,0]$ for the walk \lq started at infinity'.
 The above limits exist and constitute a probability distribution on $(-\infty,0]$ \cite[Theorem 30.1]{S}.
 It holds  that $H_{(-\infty, 0]}^{+\infty}(y) $ is bounded above and below by positive multiples of  $\int_{-\infty}^y F(u)du$; in particular 
 \beqn\label{if_if}
\sum_{y=-\infty}^{-1} H_{(-\infty, 0]}^{+\infty}(y)(-y) <\infty \;\; \mbox{ if and only if } \;\; 
  E[|S_1|^3; S_1<0]<\infty
  \eeqn
  (cf. e.g., \cite[(2.7)]{U1dm.f.s}).
It is observed in \cite{U1dm.f.s}  that $S^{b_N}$ conditioned on event (\ref{CD}) either comes  near to $A$ but still avoids it or clears $A$ by one   long jump that becomes indefinitely large as $N\to\infty$  according as  the above infinite series  is convergent or divergent.
This would convince one  that if the series are convergent there appears a continuous process in the limit.
In case of the  divergence  the limit process  may still be continuous and in order to ascertain  it to be discontinuous   we need to estimate the length of a typical   jump to reveal  whether it  is comparable to the scale  $\sqrt N$ so as to remain positive in the limit. (An answer to the question will be given in  Lemma \ref{lem3.3}.)

\v2\v2

{\bf Results  for    a random walk with $\sigma^2=\infty$}.   
\v2

Let $1<\alpha<2$
and   the distribution  function  of $S_1$  satisfy
\beqn\label{stbl0}
1-F(u) \sim L(u)/u^\alpha \quad\mbox{and}\quad F(-u)/[1-F(u)] \to 0 \quad  \mbox{as}\quad u\to\infty,
\eeqn
where $L(u)$ is a continuous positive  function on $u\geq 0$ slowly varying at infinity. 
This condition  is necessary and sufficient in order for   the scaled process  $S_{\lfloor nt \rfloor}/\la_n$ to converge in law to   a strictly stable process of index $\alpha$ which has  no negative jump, provided that the norming constants $\la_n$ are suitably chosen,  which we may and do take so that $\la_n^\alpha /n L(\la_n)\to 1$.  Let $a(x)$ denotes the potential function of the walk $S$ (defined in (\ref{D_a})). 
Then, corresponding to the equivalence relation   (\ref{if_if}) it  holds  that  under the condition (\ref{stbl0})
\beqn\label{H/aF}
\sum_{y=-\infty}^{-1} H_{(-\infty, 0]}^{+\infty}(y)a(y) <\infty\;\;\mbox{  if and only if}\;\;
\sum_{y=-\infty}^{-1}  [a(y)]^2 F(y)<\infty.
 \eeqn
  Let  $C^+$ stand for the first sum in (\ref{H/aF}) so that $C^+<\infty \Longleftrightarrow \sum_{y=-\infty}^{-1}  [a(y)]^2F(y)<\infty$. 
Then
\beqn\label{eqL3.3}
\lim_{x\to\infty} a(x) = C^+\;  (\leq \infty).
\eeqn
 (Cf. \cite[Theorem 2 and Corollary 2]{Uladd}  for (\ref{H/aF}), (\ref{eqL3.3}).) 
Denote by $(Y_t)_{t\geq 0}$ the limiting stable process started from zero  and put $Y^x_t =x +Y_t$.
Let $\mathfrak{p}_t(x)$  be the 
probability density  of $Y_t$  and $\mathfrak{p}^0_t(x,\cdot) $ the transition  density  of the   stable process killed  on hitting  zero. For $x<0$   define 
$ \rho_t(x)  $ to be the density of the hitting time of $x$ for $Y$.  Then
$$\rho_t(x) =  \frac{x}{t}\mathfrak{p}_t(-x)\quad \mbox{for}\;\; x > 0$$
(cf. \cite[Corollary VII.3]{Bt}).  It follows that for $x<0$,
\beqn\label{ent_l}
\rho_{s+t}(-x) = \int_{-\infty}^0 \rho_t(-z) \mathfrak{p}^0_s(z, x)dz
\eeqn
(cf. Appendix (B)), saying that  the family $(\rho_t(-x)dx: t>0, x <0)$ constitutes an entrance law for the transition semigroup  of  the processes $(Y^{x}_t)_{t\geq0}, x<0$ killed on exiting from the negative half line $(-\infty,0)$. Note that the transition function of this killed process agrees with $\mathfrak{p}^0_t(x,y)$ restricted to $x\vee y<0$.

Let  $(X^{b,-c,T}_t)_{0\leq t\leq T}$  be the process whose transition law is given as before by (\ref{m1}) and (\ref{m2})  but with $\mathfrak{p}_t^0$ in place of $\mathfrak{g}_t^0$ and similarly for  $\breve X^{b,c, T}_t$, $t\leq T$. 
 Instead of (\ref{b/c}) let 
$b_N\sim b\la_N$ and $c_N\sim c \la_N$ and define the scaled process $ X^{(N)}_t$ by
\beqn\label{XX}
 X^{(N)}_t=  S^{b_N} _{\lfloor t N \rfloor} /\la_N,\quad 0\leq t\leq T.
\eeqn  
 The probability of event (\ref{CD}) is supposed to be positive for  $N$ large enough.  Then the following theorems hold for  $ X^{(N)}_t$.
It is noted  that under (\ref{stbl0}) 
 $$a(-y) \sim \kappa_\a \, y^{\alpha-1}/L(y) \quad \mbox{as}\quad y\to +\infty$$
where $\kappa_\a = 1/\Ga(\a)\Ga(\a-1)$ (cf. \cite[Lemma 3.1, Eq(9.2)]{Uattrc}, \cite[Lemma 3.3]{B2}), so that $C^+<\infty$   if $F(-y) = O(y^{-\beta})$ with  $\beta< 2\alpha-1$.

\begin{Thm}\label{thm3} \, Suppose (\ref{stbl0}) to hold and  that  either $C^+ <\infty$ or  $F(u)$ is regularly varying with index $-2\a+1$ as $u\to -\infty$.  Then the  law of $( X^{(N)}_t: 0\leq t\leq T) $ conditioned  on   event (\ref{CD})
  converges to the law of $(X^{b,-c,T}_t)_{0\leq t\leq T}$ w.r.t.  the Skorohod  topology of $D([0,T],\R)$.
\end{Thm}
\begin{Thm} \label{thm4} \, Suppose (\ref{stbl0}) to hold and that    $F(u)$ is regularly varying as $u\to -\infty$ with index $-\beta$,  $\a\leq \beta< 2\a-1$.  
     Then the   law of $(X^{(N)}_t: 0\leq t\leq T) $ conditioned  on   event (\ref{CD})  converges to the law of $(\breve X^{b,-c,T}_t)_{0\leq t\leq T}$
 w.r.t.  the Skorohod topology of $D([0,T], \R)$.
\end{Thm}

{\sc Remark 3.}  In Theorems \ref{thm1} through \ref{thm4}  the scaled walk  $X^{(N)}_t$ is considered under  the conditional law given the two events given in  (\ref{CD}). By the same token as what is mentioned in Remark 1(c) we can replace the first 
of it by $\sigma_A^{b_N}=\infty$  with the conclusion kept unaltered;  in other
words, if  $P^{(A,b)}$  denotes the conditional law $P[\cdots | \sigma_A^{b_N} =\infty] $, the bridge
$P^{(A,b)}[ \cdots | X^{(N)}_T =- c_N] $ converges to  the same  limit laws as specified above.
\v2

In the next section we collect the known  results  concerning the walk $S^x$ 
killed on  $A$  which  are  used in the proofs of Theorems \ref{thm1} to \ref{thm4},  and prove some related lemmas, especially Lemma \ref{lem3.3} mentioned above. The proofs of  Theorems \ref{thm1} and \ref{thm2} are given in Sections 4.1 and 4.2, in which the tightness of the sequence of conditional laws of  $X^{(N)}$   is verified and  the finite dimensional distributions of limit processes are derived.  Theorems \ref{thm3} and \ref{thm4}
are proved in almost the same way as Theorems \ref{thm1} and \ref{thm2} and we do  not
provide full proofs of them  except for some remarks and key lemmas that we give in Section 4.3 as well as Section 3. 

\section{Preliminary Results }

We present known results taken mainly from \cite{S}, \cite{U1dm} and \cite{U1dm.f.s} and prove some related results. In the rest of this paper  the letters  $x, y, z$ and $w$  always  denote the  integers representing states of the walk. 

For a non-empty  $B \subset \Z$   denote  by $p_B^n(x,y)$  the transition probability of the walk $S_n$ killed upon entering $B$, which we define by
 $$p^n_B(x,y)  =P[ S^x_n = y, \sigma_B^x >n], \qquad n=0,1, 2,\ldots. $$
This entails   $p_B^n(x,y) =0$ whenever $y\in B, n\geq 1$; and  $p_B^0(x,y)= \de_{x,y}$ even if $y\in B$,
where $\de_{x,y}$ equals unity if $x=y$  and zero  if $x\neq y$ 
(recall that $\sigma^x_B\geq 1$ for every $x$). Put
 $$p(x) = P[S_1=x]\quad \mbox{ and}\quad  F(u) =  P[ S_1 \leq u]  \quad (u\in \R).$$
We  suppose that the walk is (temporally) {\it aperiodic}, namely for every $x\in \Z$ there exists $n_x\geq 1$ 
such that $P[S_n =x] >0$ for $n\geq n_x$, which  does not give rise to any loss of generality.

Let $a(x)$ be 
the potential function  of the walk defined by
\beqn\label{D_a}
a(x) = \lim_{n\to\infty}\sum_{k=0}^n [p^k(0) - p^k(-x)].
\eeqn
 It holds    that 
$a(x+1) -a(x) \to  \pm 1/\sigma^2$ as $ x\to \pm\infty$,
which  implies   
$$a(x+z) -a(x) = \pm z\{1/\sigma^2+o(1)\} $$
with  $o(1) \to 0$ as $x\to \pm \infty$ uniformly for  $z$ with  $(x+z)x>0$,  and
 $a(x)/| x| \to 1/\sigma^2$
as $|x|\to  \infty$. 
(Cf.  Theorem 28.1 and Theorem 29.2 of \cite{S}.) Let $g_{\{0\}}(x,y) = \sum_{n=0}^\infty p^n_{\{0\}}(x,y)$,  the Green function of the  walk killed on visiting  the origin.  
We have the identity  
$$g_{\{0\}}(x,y)  =a^\dagger (x)+a(-y) -a(x-y)$$
where $a^{\dagger}(x) =a(x)+ \delta_{0,x}$ \cite[Proposition 29.4]{S}.
 In what follows these relations will be  used frequently and  not   noticed of their use.  On putting
$$\la(x) = a(x)- x/\sigma^2,\, x\in \Z$$
it also holds \cite[(2.9)]{U1dm}, \cite[(2.15)]{U1dm.f.s} that if $\sigma^2<\infty$,
\beqn\label{a_r}
E[ \la(S^x_{\sigma_{(-\infty,0]} })] = \la(x)\quad \mbox{for}\quad x\geq 1.
\eeqn

Let $A$ be a non-empty finite subset of $\Z$.
For convenience of  description  we assume    that 
$\max A=0$
(namely  $0\in A$ and $A\subset (-\infty,0]$)   and that
$$P[ \sigma^x_{\{y\}} <\sigma_A^x]>0 \;\;\mbox{    for some \;$x> 0$\, and\, $y< \min A$.} $$

In the following statements a constant $M>1$ is given  arbitrarily in advance.  In the square bracket at the head of each of them  is indicated  the proposition   which the  result  is taken from. $C, C', C_1, \ldots$ etc. denote unimportant positive 
 constants  whose  values  may depend on $F$   and  vary  at different occurrences of them.  Sometimes these constant depends on $M$ and possibly on $A$, in which case we write $C_M$, 
$C_{M,A}$, etc. 

\v2{\bf P0} \; \cite[Remark 4 (Sect. 5), Theorem C (Sect. 4.2) ]{U1dm.f.s} \; {\it  Uniformly  for positive  $ x, y < M\sqrt n$,} 
$$p^n_{(-\infty,0]}(x,y)  = p^n_{\{0\}}(x,y)\{1+ O(y^{-1} + \la(x)x^{-1})\} \quad \mbox{and} \quad p^n_{(-\infty,0]}(x,y)  \asymp xy/n^{3/2}.$$
(Here and in the sequel  $\asymp$ means that the ratio of two sides of it is bounded away from zero and infinity.)

In view of  the first  relation above  $p^n_{\{0\}}(x,y)$ and  $p^n_{A}(x,y) $ may be interchangeable in most of  the arguments made later if $xy>0$, since the precision of  estimates for small values of $x, y$ are irrelevant for the proofs of Theorems \ref{thm1} and \ref{thm2}.

\v2
{\bf P1} \,\cite[Theorem 1.1]{U1dm}.\,  
{\it 
{\rm (a)}  \, Uniformly for $x\in \Z$ and $y\in \Z$ subject to  the constraints  $-M\leq x\leq M\sqrt  n$  and $-M\leq  y\leq M\sqrt  n$, as  $n\to\infty$ and $(|x|\wedge |y|)/\sqrt n \to 0$
\beqn\label{eq_thmA}
p_{\{0\}}^n(x,y) = \frac{\sigma^4 a^\dagger(x)a(-y) +xy}{\sigma^2 n}\mathfrak{g}_{\sigma^2 n}(y-x)\{1+o(1)\}.
\eeqn

\v2
{\rm (b)}  \,  As $x\wedge y \wedge n \to\infty$ under $x\vee y<M\sqrt n$ 
\beqn\label{eq_thm10}
p_{\{0\}}^n(x,y)= \mathfrak{g}^0_{\sigma^2 n}(x,y)\{1+o(1)\}.
\eeqn
\v2
{\it 
{\rm (c)} \, Whenever $1\leq |x| \wedge |y|<  \frac12(|x|\vee |y|)$,  \,}
$$ p_{\{0\}}^n(x,y)    \leq C_M\frac{\sqrt n\wedge |x|\wedge |y| }{(|x| \vee|y|)^2}.$$
[Here Spitzer's  bound $p^n(x) = o(\sqrt{n}/x^2)$ ($|x|/\sqrt n\to\infty$) (cf. \cite{S}) is also employed.]
\v2
{\rm (d)}  \, For all  $n\neq 0$, $x\neq 0$ and $y$, \;\; 
$$ p_{\{0\}}^n(x,y)  \leq C|xy|/n^{3/2};$$
 in case  $(|x|\vee |y|)/\sqrt n \to\infty$ the right side may be replaced by $o(xy/n^{3/2})$. }

\v2

{\bf P2} \,\cite[Theorem 2]{U1dm.f.s}.\, {\it If  $E[ \,|S_1|^{3}; S_1<0]<\infty$, then 
\v2

{\rm (a)} \,\;  $E[\la(S^x_{\sigma_{(-\infty,0]} })]= \la(x) 
\,  \longrightarrow \, \sum^{\infty}_{z=0} H^{+\infty}_{(-\infty,0]}(-z)[a(-z)+z/\sigma^2]<\infty$ \quad as\quad $x\to\infty;$\,
\v2
{\rm (a$'$)} \, $H^{+\infty}_{(-\infty,0]}(-z) =P[\hat Z<-z]/E|\hat Z| = o(1/z^2);$ \, and

\v2
{\rm (b)} \;  as $x\wedge (-y) \wedge n\to\infty$ subject to the condition  $x\vee(- y) <M\sqrt n$
$$p_A^n(x,y) \sim C_A^+ \rho_{\sigma^2 n}(x+|y|),$$
where $C_A^+$ is some positive constant determined by  $A$ and $F$ (see \cite{U1dm.f.s} for an explicit form).} 
\v2
  
 {\bf P3} \cite[Proposition  8]{U1dm.f.s}.\, {\it For
  $-M\sqrt n <y<0< x<M\sqrt n$ and for $n$ large enough
\beqn\label{eq_pr7}
 p^n_A(x,y)  \geq c_M \bigg(1+ \sum_{w=1}^{x\wedge|y|} p(-w)w^3\bigg)\frac{x + |y|}{ n^{3/2}},
 \eeqn
where  $c_M$  is a positive constant   depending  on $M$  as well as  $F$ and $ A$.}

\v2
{\bf P4} \,\cite[Corollary 7]{U1dm.f.s}.\, 
{\it Under   $-M\sqrt n <y<0< x< M\sqrt n$   }
\beq
&& P\big[S^x_{\sigma_{(-\infty,0]}}<-\eta \,\big|\, n<\sigma^x_{\{0\}}, S^x_n =y\big] 
\\[2mm]
&&\quad  \longrightarrow  \left\{ 
\begin{array}{ll}
0 \quad\mbox{ as}  \;\;  \eta \to\infty\;\;\;\quad  \mbox{uniformly for} \;\; x, \,y\; &\mbox{ if}\quad E[|X|^3; X<0] <\infty, \\[1mm]
1 \quad\mbox{ as}  \; \;  x\wedge (-y) \to\infty \;\; \mbox{ for each}\; \eta>0\;\; &\mbox{ if}\quad E[|X|^3; X<0] =\infty.
\end{array} \right.
\eeq

  \v2
 
{\bf P5}  \, \cite[Corollary 3]{U1dm.f.s}, \cite[Propositions 2.2 and 2.3]{U1dm}.  {\it Uniformly for $x>0$ and $r>0$, as $R\to\infty$}
$$P[S^x_{\sigma_{[R,\infty)}} > R+ r; \, \sigma^x_{[R,\infty)} < \sigma^x_{(-\infty,0]} ]
\leq r^{-1}x \times o(1).$$
\v2

{\bf P6} \,\cite[Theorem 1.4]{U1dm}.  {\it  Let  $h_x$  be the space-time distribution of the first entrance of $S^x$ into $(-\infty,0]$, namely
 $h_x(n,y) =P[S^x_{\sigma_{(-\infty, 0]}}=y, \sigma^x_{(-\infty, 0]}=n]  \;\; (x\geq 1, y\leq 0).$
Then  uniformly   for $y\leq 0<x < M\sqrt n $, as $n\to\infty$
\beqn\label{h0}
h_x(n,y)= \frac{ f^+(x)\mathfrak{g}_n(x)}{n} H_{(-\infty, 0]}^{+\infty}(y)\{1+o(1)\} + o\bigg(\frac{x}{n^2}\bigg);
\eeqn
and  there exists a positive constant  $c$ such that whenever  $y\leq 0, x\wedge n\geq 1$,
\beqn\label{h00}
h_x(n,y) \leq c H_{(-\infty, 0]}^{+\infty}(y)\bigg(\frac1{n}\wedge  \frac1{x\sqrt n} \bigg).
\eeqn
Here  $H_{(-\infty, 0]}^{+\infty}$ is the probability law given in (\ref{H+})  and $f^+(x),\, x\geq 1$ is a positive multiple of the renewal function of the descending ladder height process with the multiplicative constant chosen so  that $f^+(x) \sim x$ as $x\to\infty$.
} [Although  (\ref{h00}) is restricted to $x>\sqrt n$ in \cite{U1dm}, the extension to $x\leq \sqrt n$ is easily verified (cf. \cite[Lemma 4.1]{U1dm.f.s} or proof of \cite[Lemma 6.5]{Uattrc}.) ]

\v2

\begin{lem}\label{lem3.1} \,  If $E[|S_1|^3; S_1<0] <\infty$, then for some constant $C$,
\[
p_{\{0\}}^n(x,y) \leq C\frac1{y^2\vee n}\qquad \mbox{whenever}\;\; y<0\leq x.
\]
\end{lem}
\v2\n
\pf \,  First consider the case  $y <- 2\sqrt n$ and  suppose that $n$ is even for convenience of description.  In the decomposition
\beqn\label{eqL3.1} 
p_{\{0\}}^n(x,y)= \sum_{k=1}^n\sum_{z=1}^\infty h_x(k,-z) p_{\{0\}}^{n-k}(-z,y)
\eeqn
we break the  double sum on the right into three parts:
\[
I := \sum_{k=1}^n\sum_{z\leq |y|/2},\quad J_{\leq n/2}:=    \sum_{k=1}^{n/2}\sum_{z > |y|/2}\quad \mbox{and}\quad J_{>n/2} := \sum_{k=1+n/2}^n\sum_{z > |y|/2} \quad \mbox{(say)}.
\]
By {\bf P1}(c)  it follows that  for $z\leq |y|/2$, $p_{\{0\}}^{n-k}(-z,y) \leq C(z\wedge \sqrt {n-k})/ y^2 \leq C' z/y^2$, so that uniformly for $x>0$,
\[ I\leq C'\sum_{0< z\leq |y|/2} H^x_{(-\infty,0]}(-z) \frac{z}{y^2} \leq C''/y^2,\]
where for the last inequality we use the assumption $E[|S_1|^3; S_1<0] <\infty$ (see {\bf P2}(a)).  Observe that  
$ \sup_{k\leq n/2} p_{\{0\}}^{n-k}(-z,y)$ is dominated by a constant multiple of  $[\sqrt n /(y+z)^{2}]\wedge n^{-1/2}$, whose sum over $z>|y|/2$ is bounded by a constant.  Using this and {\bf P2}(a$'$) in turn we deduce that
\[ J_{\leq n/2} \leq \sum_{z>|y|/2}H^x_{(-\infty,0]}(-z)  \sup_{k\leq n/2} p_{\{0\}}^{n-k}(-z,y)\leq C\sup_{z>|y|/2} H^{+\infty}_{(-\infty,0]}(-z) = o(1/y^2).
\]
By the second half of {\bf P6}  we have $h_x(k,-z) \leq Cn^{-1}H^{+\infty}_{(-\infty,0]}(-z) $ for $k>n/2$, and hence, employing  {\bf P2}(a$'$) again, 
\[ J_{> n/2} \leq \frac{o(1/y^2)}{n}\sum_{k=n/2}^n \sum_{z>|y|/2} p_{\{0\}}^{n-k}(-z,y) =o(1/y^2).
\]
Thus the bound of the lemma is obtained if $|y|\geq 2\sqrt n$. 

For $0\leq y< 2\sqrt n$, we break the outer sum at $n/2$  in the right side expression of  (\ref{eqL3.1}). Use  (\ref{h00}) together with $g_{\{0\}}(-z,y) \leq g_{\{0\}}(-z,-z) \leq Ca(-z)$ ($z>0$)   for the sum over $k > n/2$,  which is then evaluated to be at most   a constant multiple of $n^{-1}\sum_{z=1}^\infty H^\infty_{(-\infty,0]}(-z)z \leq C/n$. The same bound is obtained for the   other sum  by  using 
  $p^{n-k}(-z,y) \leq C zy/ n^{3/2}$ valid for $k\leq n/2$  (cf. {\bf P1}(d)).  \qed

\begin{lem}\label{lem3.3} \,Suppose that $F$ is regularly varying at $-\infty$ with index  $-\beta$, $\beta\in [2,3]$. Then for each $M > 1$  the following hold uniformly for $x, y$ satisfying  $-M\sqrt n  <y< 0<x <M\sqrt n $.

\v2
 {\rm (a)}\;\;If $2\leq \beta <3$,  
  $$   P\big[S^x_{\sigma_{(-\infty,0]}}\leq -\e (x \wedge |y|)\, \big| \; S^x_n =y, \sigma^x_A> n \big] \,\to \,1 \quad \mbox{as}\;\; x\wedge (-y) \wedge \e^{-1}\to\infty.$$ 
 
 {\rm (b1)} \,If $2<\beta <3$,
 \beqn\label{eq_b1}
P\big[S^x_{\sigma_{(-\infty,0]}} < - K (x\wedge |y|)\, \big| \; S^x_n =y, \sigma^x_A> n \big]\, \,\to \,0 \quad   \mbox{as} \;\; x\wedge (-y) \wedge K \to\infty.
\eeqn

 {\rm (b2)} \,If $\beta=2$, 
 (\ref{eq_b1}) holds for $x,y$ satisfying   $ x\wedge |y| \geq M^{-1}\sqrt n$ in addition.
 
 \v2
 {\rm (c)} \;\;If $\beta=3$ and  $E\big[|S_1|^3;S_1<0\big]=\infty$, then for each $\e>0$ 
 $$  
P\big[  -  \e (x \wedge |y|) < S^x_{\sigma_{(-\infty,0]}} <- 1/\e \; \big| \; S^x_n =y, \sigma^x_A> n \big] \,\to\,1\quad \mbox{as} \;\; x\wedge |y|\to\infty.
$$
\end{lem}

\v2\n
\pf\,  Let $F$ be  regularly varying as is assumed in the lemma. As $m\to\infty$ by Karamata's theorem
 \beqn\label{Kar}
  \sum_{w=1}^m p(-w)w^3 \sim \frac{\beta}{3 - \beta}m^3F(-m)\quad (\beta\neq 3).
  \eeqn
 Note that $P[S^x_{\sigma_{(-\infty,0]}}\geq -\e (x\wedge|y|), \, S^x_n =y, \sigma^x_A>n ]$ is non-increasing  in  $A$  and then observe    that (a) follows if we  show that  as $x\wedge (-y)\wedge n\to\infty$ 
\beqn\label{eq1}
\frac{\sum_{k=1}^n\sum_{w=1}^\infty \sum_{-\e (x\wedge |y|) \leq z<0} p^{k-1}_{(-\infty,0]}(x,w)p(z-w)p^{n-k}_{\{0\}}(z,y)}
{(x\vee|y|)n^{-3/2}   \sum_{w=1}^{x\wedge |y|} p(-w)w^3 } \,\longrightarrow\, 0, 
\eeqn
owing to {\bf P3}.
As before we split the outer summation of the numerator at $k=n/2$. 
Employing the  obvious inequality $p^{k-1}_{(-\infty,0]} \leq p_{\{0\}}^{k-1}$ together with  {\bf P1}(d)  we infer that 
the sums over  $k\leq n/2$ and $k>n/2$ are dominated by constant multiples of
\beqn\label{eq3}
|y|n^{-3/2} \sum_{-\e (x\wedge |y|)\leq z<0}\; \sum_{w=1}^\infty g_{\{0\}}(x,w)p(z-w) |z| 
\eeqn
and
\beqn\label{eq4}
xn^{-3/2} \sum_{-\e (x\wedge |y|)\leq z<0} \;\sum_{w=1}^\infty  w p(z-w) g_{\{0\}}(z,y), 
\eeqn
respectively.  We have  the bound 
 $g_{\{0\}}(x,y) \leq C (|x|\wedge |y|)$ for all \;$x$\; and\; $y$\;  with\; $xy\neq 0$.
 On putting  $m=\lfloor \e (x\wedge |y|)\rfloor$ and  substituting from  $g_{\{0\}}(x,w)\leq C|w|$ and  $g_{\{0\}}(z,y)\leq C|z|$  an elementary   computation verifies that both of the sums of the above double series are dominated by a constant multiple of $\sum_{w=1}^\infty\sum_{z=1}^m p(-z-w)zw$ which is at most
\beqn\label{eq5}
\frac12 \sum_{u=1}^{2 m}\sum_{v=0}^u p(-u)(u^2-v^2) + \frac12  m^2\sum_{ m}^\infty p(-w)w \leq C m^3F(- m).
\eeqn
Since  $ m^3F(- m) \sim \e^{3-\beta}(x\wedge |y|)^{3}F(-x\wedge|y|)$,  this together with (\ref{Kar})   shows  (\ref{eq1}) if $\beta<3$. 

As for (b1) and (c)  we show shortly  that if $\beta>2$,  then for all  $m$ large enough
\beqn\label{eq2}
\sum_{k=1}^n \sum_{z <- m} p^{k-1}_{\{0\}}(x,w)p(z-w)p^{n-k}_{\{0\}}(z,y)\leq \frac{C(x\wedge |y|)m^2F(-m) (x\vee |y|)}{n^{3/2}}.
\eeqn
If $\beta\neq 3$, on taking  $m=\lfloor K(x\wedge |y|)\rfloor$   (b1) is immediate from this  in view of (\ref{Kar}).  
If $\beta=3$, then $w^3F(-w)$  is slowly varying  and    it follows  that  $\sum_1^m w^3 p(-w) \sim  3\sum_1^m w^2F(-w)$ and   
$$\frac{m^3F(-m)}{\sum_1^m w^3 p(-w)} \to 0\quad \mbox{ as}\quad  m\to \infty.$$ 
On using {\bf P3} as before this  together  with  (\ref{eq2}) shows that the conditional probability of the event $S^x_{\sigma_{(-\infty,0]}} \leq -\e (x\wedge |y|)$ tends to zero.  (c) now follows from {\bf P4}. 

 The proof of (\ref{eq2}) is similar to the  one given  above for  (\ref{eq1}).  What we should   evaluate  are 
\beqn\label{eq6}
\begin{array}{ll}
 &\sum_{z< -  m}\; \sum_{w=1}^\infty g_{\{0\}}(x,w)p(z-w)p_{\{0\}}^n(z,y); \;\mbox{and}\\[4mm]
 &\sum_{z< -  m}\; \sum_{w=1}^\infty p_{\{0\}}^n(x,w)  p(z-w)g_{\{0\}}(z,y)
 \end{array}
\eeqn 
instead of  (\ref{eq3}) and  (\ref{eq4}), respectively. By   $g_{\{0\}}(x,w) \leq Cx$ it plainly follows
 that
\beqn\label{eq_c}
 \sum_{w=1}^\infty g_{\{0\}}(x,w)p(z-w) \leq C' x F(z),
 \eeqn
and  if $\beta>2$    the sum of the first double series  in (\ref{eq6}) is at most a constant multiple of
\beqn\label{eq_c1} x\sum_{z< - m}F(z)p_{\{0\}}^{n}(z,y)\leq C x|y| m^2F(- m)n^{-3/2}. 
\eeqn
 Noting that $\sum_{w=1}^\infty wp(z -w) =\sum_{w=1}^\infty F(z-w)\sim (\a-1)^{-1}F(z)|z|$ as $z\to-\infty$  a similar computation leads to the same  upper bound for the second one. This shows (\ref{eq2}),  since  $|xy| = (x\wedge |y|)(x\vee |y|)$. 

It remains to verify (b2). If $\beta=2$, the upper bound (\ref{eq_c1}) is not valid but in (b2)  it is supposed that $x\wedge |y|>\sqrt n/M$ so that if  $m=K(x\wedge |y|)$ ($K>1$), then    $p_{\{0\}}^{n}(z,y) \leq C\sqrt n /z^2$ for $z<-  m$  and $|y|<M\sqrt n$  owing to  {\bf P1}(c), and we obtain  instead of
 (\ref{eq_c1}) 
 $$x\sum_{z< - m}F(z)p_{\{0\}}^{n}(z,y)\leq C' x\sqrt n \sum_{z<-m} F(z)/z^2\leq \frac{C_M x F(-m)}{K}\leq  C_M (x\vee|y|) F(-m).
$$
Observing  that the second sum in (\ref{eq6}) admits the same upper bound we   conclude
 (b2) to be  true since $F(-m)\leq 
C_M'K^{- 2} (x\wedge|y|)^3 F(- x\wedge |y|)/n^{3/2}$  for $n$  large enough. 
\qed

\begin{Cor}\label{cor3.1}\, Let $F$ be as in the preceding lemma and suppose $E[|S_1|^3;S_1<0]=\infty$ if $\beta =3$.  Then for each $M > 1$,  uniformly for $x, y$ satisfying  $-M\sqrt n  <y< 0<x <M\sqrt n $, 
\v2
{\rm (a)}\quad 
$P\big[ S^x_k <0 \;\; \mbox{for}\;\; \sigma^x_{(-\infty,0]} \leq k \leq n   \,\big|\, S^x_n=y, \sigma^x_A> n\big]  \; \to \; 1$ \qquad $(x\wedge|y|\to\infty)$,
\v2
{\rm (b)}\;   if\; $2\leq \beta<3$, 
$$\qquad
P\big[ \, |S^x_{\sigma_{(-\infty,0]}}| \wedge S^x_{\sigma_{(-\infty,0]}-1} >  \e (x\wedge|y|)  \,\big|\, S^x_n=y, \sigma^x_A> n\big]  \to 1\qquad (x\wedge |y| \wedge \e^{-1} \to\infty),
$$

\v2
{\rm (c)}\;  if\; $\beta=3$, then for each $\e>0$ 
$$\quad
P\big[ \, |S^x_{\sigma_{(-\infty,0]}}|\vee S^x_{\sigma_{(-\infty,0]}-1} < \e (x\wedge|y|)  \,\big|\, S^x_n=y, \sigma^x_A> n\big] \; \to \; 1  \qquad (x\wedge |y|\to\infty).
$$
\end{Cor} 
\v2\n
\pf\, The first relation  follows from ${\bf P4}$ and ${\bf P0}$, the latter asserting $p^n_{\{0\}}(z,y)/ p^n_{[0,\infty)}(z,y) \to 1$
as $z\vee y \to -\infty$. By duality relations the rest follows from Lemma \ref{lem3.3}. \qed 
\v2

{\bf Corresponding  lemmas  for the case $\sigma^2=\infty$}
\v2

Here we suppose that (\ref{stbl0}) is satisfied and   prove   Lemmas \ref{lem3.4} and \ref{lem3.5}  below  that correspond to  Lemmas \ref{lem3.1} and  \ref{lem3.3}, respectively,  and will be  used   for the proofs of Theorems \ref{thm3} and \ref{thm4}.  We shall use the following large deviation estimate:
\beqn\label{LD}
  p^n(x) \leq C\{1\wedge [nL(|x|)/{|x|^\alpha}]\}/\la_n 
  \eeqn
 (cf.  \cite{Brg}).  
 We also have 
\beqn\label{H/aa}
H^x_{(-\infty,0]}\{a\} =a(x), \quad  x\geq 1
\eeqn
  valid if $E [ S^0_{\sigma_[0,\infty)}]=\infty$ always true under   (\ref{stbl0})  (cf. \cite[Corollary 1]{Uladd}). Recall $\lim_{x\to\infty} a(x) = C^+ \leq \infty$. Most of the results given below are based on \cite{Uattrc} in which the transition function of the killed walk denoted by $Q_A^n(x,y)$  is  defined slightly differently from but  agrees with  $p^n_A(x,y)$ whenever  $y\notin A$. 
 
\begin{lem}\label{lem3.4} \,  If $C^+ <\infty$, then for some constant $C$,
\[
p_{\{0\}}^n(x,y) \leq C\bigg(\frac{L(|y|)}{(-y)^\alpha}\wedge \frac1{n}\bigg)\qquad \mbox{whenever}\;\; y<0\leq x.
\]
\end{lem}
\v2\n
\pf \,   
For $y<-2\la_n$, the proof parallels to that of Lemma \ref{lem3.1} with the help of the following bounds 
\beqn\label{p_p/y}
p^n_{\{0\}}(x,y)\vee p^n_{\{0\}}(y,x)\leq C\frac{[\,|x|^{\alpha-1}/L(|x|)]\wedge [n/\la_n]}{ |y|^\alpha/L(|y|)} \quad \; \;  (1\leq |x|< |y|/2),
\eeqn
  \beqn\label{h/H}
  h_x(n,y) \leq C\bigg(\frac1{n}\wedge\frac{L(x)}{x^\alpha} \bigg)\frac{x}{\la_n}H^{+\infty}_{(-\infty]}(y)    \quad \;\;  ( y<0 <x)
  \eeqn
that follow immediately  from Lemma 6.2 (combined with (\ref{LD}))  and  Lemma 6.5(ii) of \cite{Uattrc},  respectively.  As  for  $-2\la_n \leq y<0$, in view of the identity 
  $p^n_{\{0\}}(x,y)=p^n_{\{0\}}(-y,-x)$ what has been just proved yields the bound  $p^n_{\{0\}}(x,y) <C/n$ 
  for $x>2\la_n$.  For $x \leq 2\la_n$, a better bound is obtained in \cite[Proposition 2.3(i)]{Uattrc}
  (applied to $p^n_{\{0\}}(-y,-x)$). 
\qed

\begin{lem}\label{lem3.5} \, Suppose  $F$ is regularly varying at $-\infty$ with index  $-\beta$, $\a\leq \beta \leq 2\a-1$ in addition to  (\ref{stbl0}). For each $M > 1$  the following holds uniformly for $x, y$ satisfying 
 $-M\la_n  <y< 0<x <M\la_n$.

\v2
 {\rm (a$^*$)}\;\;If $\a\leq \beta <2\a-1$,  
 $$   P\big[S^x_{\sigma_{(-\infty,0]}}\leq -\e (x \wedge |y|)\, \big| \; S^x_n =y, \sigma^x_A> n \big] \,\to \,1 \quad \mbox{as}\;\; x\wedge (-y) \wedge \e^{-1}\to\infty.$$ 
 
 {\rm (b$^*$1)} \,If $\a <\beta <2\a-1$, under the additional constraint $x\vee |y| \leq M(x\wedge |y|)$
 \beqn\label{eq_b*1}
P\big[S^x_{\sigma_{(-\infty,0]}} < - K (x\wedge |y|)\, \big| \; S^x_n =y, \sigma^x_A> n \big]\, \,\to \,0 \quad   \mbox{as} \;\; x\wedge (-y) \wedge K \to\infty.
\eeqn

 {\rm (b$^*$2)} \,If $\beta=\a$, 
 (\ref{eq_b*1}) holds for $x,y$ satisfying   $ x\wedge |y| \geq M^{-1} \la_n $ in addition.
 
 \v2
 {\rm (c$^*$)} \;\;If $\beta= 2\a-1$ and  $C^+=\infty$, for each $\e>0$ 
 $$  
P\big[  -  \e (x \wedge |y|) < S^x_{\sigma_{(-\infty,0]}} <- 1/\e \; \big| \; S^x_n =y, \sigma^x_A> n \big] \,\to\,1\quad \mbox{as} \;\; x\wedge |y|\to\infty.$$
\end{lem}
\v2
The proof of the lemma proceeds parallel to that of Lemma \ref{lem3.3} and we point out only  main steps after    stating  the results from \cite{Uattrc}  needed for it. 
Put for $x>0$
$$D_n(x)=  \frac{a(x)\la_n}{n^{2}}  + \frac{x}{ n\la_n}.$$
Instead of {\bf P3} it holds that if  either  $C^+ <\infty$ or $F(x)$ is regularly varying as $x\to-\infty$, then 
\beqn\label{eqP1}
p^n_{\{0\}}(x,y) \geq C_M^{-1}\{D_n(x) a^\dagger(-y)+D_n(-y) a^\dagger(x)\} \;\quad (-M\la_n <y \leq 0\leq x <M\la_n)
\eeqn
(the reversed inequality with  $C_M^{-1}$ replaced by $C_M$ also holds), where $C_M$ is a positive constant (see Proposition 2.3(ii) of  \cite{Uattrc} as well as the comment given right after it). We also have
\beqn\label{eqP3}
p^n_{\{0\}}(x,w) \leq   C D_n(x)\{a(-w) \wedge a(-\la_n)\} \quad  (0<x <M\la_n, w>0)
\eeqn
and its dual $p^n_{\{0\}}(z,y)  = p^n_{\{0\}}(-y,-z)  \leq C D_n(-y)\{a(z) \wedge a(-\la_n)\}$ $(z<0, -M\la_n < y<0)$ (cf. \cite[Lemma 6.1(i)]{Uattrc}).  Since $D_n(x) \leq C[x^{\alpha-1}/L(x)] \la_n/n^2$,  (\ref{eqP3}) together with $p^n(x) \leq C/\la_n$ entails    
 \beqn\label{eqP2}
p^n_{\{0\}}(x,y) \leq C \frac{|xy|^{\alpha-1} \la_n}{L(|x|)L(|y|) n^2} \quad (xy\neq 0).
\eeqn

 Proof of (a$^*$).\; 
  (\ref{eq1}), (\ref{eq3}) and (\ref{eq4}) are modified in an obvious way according to (\ref{eqP1}),  (\ref{eqP2}) and (\ref{eqP3}). Put $m=\lfloor \varepsilon(x \wedge |y|) \rfloor$ as before.
   Noting $g_{\{0\}}(x,y) \leq C\big[a(-|x|)\wedge a(-|y|)\big]$ ($xy\neq 0$), we see first that 
   $\sum_{w=1}^\infty\sum_{z=1}^m p(-z-w)a(-z)a(-w)$ is at most a constant multiple of 
  $$
 \sum_{u=1}^{2  m}\sum_{v=0}^u \frac{p(-u)(u^2-v^2)^{\a-1}}{L(u+v)L(|u-v|)} +    \frac{m^{\a}}{L(m)}\sum_{m}^\infty \frac{p(-w)w^{\a-1}}{L(w)} \leq C\frac{ m^{2\a-1}F(- m)}{L^2(m)},
$$
an estimate corresponding  to  (\ref{eq5}), and then  that (a$^*$) follows if we show that
  \beqn\label{pra*}
  \frac{a(x)\vee a(-y)}{n^2/\la_n}\cdot \frac{ m^{2\a-1}F(- m)}{L^2(m)} \bigg/
  \frac{a(x)a(-y)}{n^{2}/\la_n} \to 0
  \eeqn
  as $x\wedge (-y)\wedge \e^{-1}\to\infty$.
Now  let $F(-u)\sim L_-(u)/u^\beta$ ($u\to\infty$) for a  slowly varying function  $L_-$.  Then for $\alpha \leq \beta<2\alpha-1$ we have
$$a(x) \sim C_1 L_-(x)x^{2\alpha-\beta -1}/L^2(x)\sim  C_1 x^{2\alpha -1}F(-x)/L^2(x) \quad (x\to\infty)$$
with a constant $C_1>0$ (cf. \cite[Proposition 6.2, Eq(6.17)]{Upot}). Recalling $m =\e (x\wedge |y|)$ we find  (\ref{pra*}) to be true, provided that $\alpha\leq \beta<2\alpha -1$.

Proof of  (b$^*$1) and (c$^*$). \,  Instead of (\ref{eq2})  one shows  that if $\beta>\a$,  for   $m \leq M\la_n$ large enough, 
\beqn\label{eq20}
\sum_{k=1}^n \sum_{z <- m} p^{k-1}_{\{0\}}(x,w)p(z-w)p^{n-k}_{\{0\}}(z,y)\leq Cm^{\a}F(-m) \{a(y)D_n(x)+a(-x) D_n(-y)\}
\eeqn
 the proof being   the same as before except for a minor modification. (b$^*$1)  readily follows from (\ref{eq20}).   As for (c$^*$) we have 
the same result  as in  {\bf P4} but with  $C^+=\infty$ instead of  $E[|x|^3; X<0]=\infty$ as well as  with $\la_n$ replacing $\sqrt n$ (cf. \cite[Proposition 2.1]{Uattrc}) and combining it with  (\ref{eq20}) leads to the assertion  as before.

Proof of  (b$^*$2). \, Using the  bound $p^{n}_{\{0\}}(z,y)\leq CD_n(-y)a(-\la_n)$ valid  for  $z< - \la_n$ (see the dual of (\ref{eqP3}))  one can  proceed as before. 
 
\section{Proofs of Theorems} 

Let $X^{(N)}$ be as in (\ref{X}) and denote by
 $P^*_{(c,N)}$ the conditional probability law of the walk $S$, and hence of $X^{(N)}$, given the event  (\ref{CD}): 
$$P^*_{(c,N)}[\,  \cdots] = P[\, \cdots \,|\,  S^{b_N}_{\lfloor NT\rfloor} =-c_N, \sigma_A^{b_N}>NT].$$ 
Let  $\zeta =\zeta_N$ and  $\zeta'=\zeta_N'$ stand for the first  time  $S^{b_N}$  enters  $(-\infty,0]$ and the last time $\leq NT$ it leaves $[1,\infty)$, respectively:  
 \beqn\label{T000}
 \zeta=\sigma^{b_N}_{(-\infty, 0]} \quad\mbox{and}\quad  \zeta'= \max \{n \leq NT: S^{b_N}_n\in [1,\infty)\}.
 \eeqn
The proofs of  Theorems \ref{thm1} and \ref{thm2} are given separately according to whether $E[|S_1|^3; S_1<0]$ is finite or infinity. We continue to use the notation of the preceding section.  For simplicity we shall suppose  $b_N =\lfloor b\sqrt {\sigma^2N}\rfloor $ and  $c_N = \lfloor -c\sqrt{  \sigma^2N}\rfloor $.

\v2
\subsection{Case $E[ |S_1|^3; S_1<0]<\infty$} 
Under the conditional law $P^*_{(c,N)}$  once the walk $S^{b_N}$ goes down below the level $-\e \sqrt N$ with  $\e>0$ prescribed,  it never takes the positive  value up to the time it terminates at $-c_N$ with the probability that approaches unity as $N\to\infty$  as we shall see (the second half of Lemma \ref{lem4.1}). On taking this for granted  
the convergence of the finite dimensional distribution of $X^{(N)}$ under $P^*_{(c,N)}$
follows immediately from {\bf P1}(b) and {\bf P2}(b): Noting $p_A^n(x,y) \sim p_{\{0\}}(x,y)$ as $|x|\wedge |y| \to\infty$ by {\bf P0},  we infer that given  $\eta_k< 0< \xi_j$ ($j=0, \ldots, m$, $k=0,\ldots, m'$)  with $\xi_0= b$ and $\eta_{m'} = -c$ and   $0=s_0<s_1<\cdots<s_m<t_{0}<\cdots<  t_{m'}:=T$,  we let
 $x_j =\lfloor \xi_j \sigma\sqrt N\rfloor$,   $y_k =\lfloor \eta_k \sigma\sqrt N\rfloor$, $n_j =\lfloor N s_j\rfloor$ and $n'_k =\lfloor N t_k\rfloor$. Then, as $N\to\infty$
\begin{eqnarray}\label{T00}
&&
(\sqrt{\sigma^2 N}\,)^{m+m'}P^*_{(c,N)}\big[S^{b_N}_{ n_j} =x_j \;\; \mbox{for}\;   j \geq 1 \;  \; \mbox{and} \; \; 
S^{b_N}_{n_k}= y_k \;\; \mbox{for}\;   k \leq m'-1 \, \big]   
\nonumber\\
&&= \frac{(\sqrt{\sigma^2 N}\,)^{m+m'}}{p_A^{m+m'}(x_0, y_{m'})}\prod_{j=1}^m p_A^{n_j-n_{j-1}}(x_{j-1}, x_{j}) \cdot p_A^{n'_0 -n_m}(x_m, y_{0})   \prod_{k=1}^{m'} p_A^{n'_{k}- n'_{k-1}}(y_{k-1}, y_{k}) \nonumber \\
&&\quad  \longrightarrow\;  
\frac{\prod_{j=1}^m \mathfrak{g}^0_{s_j-s_{j-1}}(\xi_{j-1}, \xi_{j}) \cdot \rho_{t_0-s_m}(\xi_m-\eta_0) \prod_{k=1}^{m'} \mathfrak{g}^0_{t_{k}-t_{k-1}}(\eta_{k-1}, \eta_{k})}{\rho_{T}(b+c)}, 
\end{eqnarray}
of which the right side is the density of the corresponding finite dimensional distribution of the limit process
 $X^{b,c,T}$ given by (\ref{X0})---with apparent change of letters $x, y$ to $\xi, \eta$.

We need to show  the 
{\it tightness}  of the law of $X^{(N)}$ under $P^*_{(c,N)}$, or what amounts to the same thing \cite{B}, that
 for any $\e>0$ there exist $\de>0$ and $N_0$ such that if $N\geq N_0$, then
\beqn \label{T0}
P^*_{(c,N)}\bigg[ \sup_{0\leq s<t\leq T}\,  \sup_{|t-s|<\de} |X^{(N)}_t-X^{(N)}_s| > \e \bigg] <\e.
\eeqn
For $0\leq j<k\leq NT$ let $\La_{j,k}(\de)$ denote  the random variable 
$$\La_{j,k}(\de) = \sup_{\frac1{N}j \leq s<t\leq \frac1{N}k ;\, |t-s|<\de} |X^{(N)}_t-X^{(N)}_s|.$$
It suffices to show, instead of (\ref{T0}), that
$$P^*_{(c,N)}[\La_{0,\zeta}(\de) >\e] \vee P^*_{(c,N)}[\La_{\zeta,\zeta'}(\de) >2\e, \zeta<\zeta']   \vee  P^*_{(c,N)}[\La_{\zeta', NT}(\de)>\e] <\e.$$ 

For simplicity we let $T=1$.  
First we show that
\beqn\label{T1}
P^*_{(c,N)}[\La_{0,\zeta}(\de)> \e] = P^*_{(c,N)}\bigg[ \sup_{0\leq s<t\leq \zeta;\, |t-s|<\de} |X^{(N)}_t-X^{(N)}_s| >  \e \bigg] <\e.
\eeqn
The conditional probability on the left side is expressed as
$$P^*_{(c,N)}[\La_{0,\zeta}(\de) > \e] = \sum_{k=1}^N\sum_{z=1}^\infty P[\La_{0,\zeta}(\de) >\e, \zeta =k, S^{b_N}_k=-z]\frac{p_A^{N-k}(-z,-c_N)}{p_A^{N}(b_N,-c_N)}.$$
Splitting the inner summation at  $M$, by {\bf P4} we can find $M>1$ such that the contribution to the double sum from  $z> M$ is less than $\e/2$, whereas
by  {\bf P1}(c, d) it follows that for any $M>1$,  
\[ 
\max_{k\leq N}\max_{1\leq z\leq M} p_A^{N-k}(-z,-c_N) \leq CM/N.
\] 
 Since $p_A^N(b_N,-c_N)\geq C'/N$ according to {\bf P2}, the latter bound shows that the ratio under the double summation sign is bounded by a constant $C_M$ for $z\leq M$, which together with the former one 
   leads to 
$$P^*_{(c,N)}[\La_{0,\zeta}(\de)>\e] \leq C_M P[ \La_{0,\zeta}(\de) >\e] +\e/2,$$
hence by the invariance principle, which entails that  $P[\La_{0,\zeta}(\de) >\e] <\e/2C_M$ for $\de>0$ small enough, we obtain (\ref{T1}).

 The bound  $P^*_{(c,N)}[\La_{\zeta', N}(\de)>\e] <\e$ follows from  (\ref{T1}) as the dual assertion  of it.

 It remains to dispose of  $P^*_{(c,N)}[\La_{\zeta, \zeta'}(\de),  \zeta<\zeta']$.   Since $\La_{\zeta, \zeta'}(\de)$ can not exceed $2\e$   if  the path of the walk  is confined in the interval   $[-\e \sqrt N, \e \sqrt N]$ for the time duration between $\zeta$ and $\zeta'$,  it suffices to show  
 $\lim_{N\to\infty} P^*_{(c,N)}[\,  |S^{b_N}_k| \leq  \e\sqrt N \;\, \mbox{for} \;\,  k:\zeta\leq k\leq  \zeta'\,] = 1,$
 which follows from the next lemma. 

 \begin{lem}\label{lem4.1} \;  If $E[|S_1|^3; S_1<\infty]<\infty$,
 then there exists a constant $C$ depending on $b, c, A$ and $F$ such that for $R>1$
   \beqn\label{T2}
P^*_{(c,N)}\big[ \, S^{b_N}_k> R \;\; \mbox{for some} \;\;  \zeta\leq k< N\,\big] \leq C/R
\eeqn
and $P^*_{(c,N)}[ \,  S^{b_N}_k< - R \;\; \mbox{for some} \;\; k < \zeta' \,] \leq C/R$.
 \end{lem}
\v2\n
\pf\,   We have only to verify  (\ref{T2}), the other relation being    its dual. Let $E[|S_1|^3; S_1<\infty]<\infty$.   Then   $p^{N}_{A}(b_N, -c_N)\sim C_{A,b,c}/N$  
according to   {\bf P2}(b), so that the conditional  probability in (\ref{T2}) is dominated by a constant multiple of 
\[
N\sum_{k=1}^N\sum_{ z< 0}h_{b_N}(k,z) 
P\big[\sigma^{z}_{[R,\infty)}<N-k< \sigma_{\{0\}}^{z},  S^{z}_{N-k}=-c_N\big], \\
\]
 which  is obviously less than
\beqn\label{T3}
 N \sup_{x\geq R, \,0<k<N}p^{N-k}_{\{0\}}(x,-c_N)
 \sum_{ z< 0}H_{(-\infty,0]}^{b_N}(z)P\big[\sigma_{[R,\infty)}^{z}<\sigma^{z}_{\{0\}}\big].
 \eeqn
We know that   $P[\sigma_{[R,\infty)}^{z}<\sigma^{z}_{\{0\}}] \sim  P[\sigma_{\{R\}}^{z}<\sigma^{z}_{\{0\}}]$ as $R\to\infty$  uniformly in $z\in \Z$ whenever $\lim_{x\to\infty} a(z)/a(-x)=0$ (cf. \cite[Proposition 5.2(i)]{Upot}), while 
\[
  P\big[\sigma_{\{R\}}^{z}<\sigma^{z}_{\{0\}}\big] =g_{\{0\}}(z,R)/g_{\{0\}}(R,R)\leq C'a(z)/R, 
 \]
so that  by (\ref{a_r}) the sum in (\ref{T3}) is $O(1/R)$, $\la(b_N)$ being bounded under the present moment condition.
 Since 
the supremum in (\ref{T3})  is dominated by a constant multiple of $1/c_N^2\sim 1/c^2\sigma^2N$ owing to  Lemma \ref{lem3.1}, we can  conclude the asserted bound. \qed

\v2
  In the next subsection we shall need  a tightness result  with the
 conditional  of avoiding $(-\infty,0]$ instead of $A$. The following result, of which  the condition  $E[|S_1|^3; S_1<\infty]<\infty$ is irrelevant, actually  shows that the random walk bridge conditioned to stay positive weakly converges to a standard Brownian meander of length 1 pinned at a prescribed point $b$ at time 1 locally  uniformly in  $b>0$. 

\begin{lem}\label{lem4.2} For each $0<\e <1/2$, as $N\to\infty$ and $\de\downarrow 0$ independently
  $$ P\big[\, \La_{0,n}(\de) > \e\, \big|\, S_n^{b_N} =x, n<\zeta \,\big] \to 0 $$
 uniformly for  $\e N <n< N$,  $0< x < \e^{-1}\sqrt N$ and $\e<b < \e^{-1}$.
\end{lem}
\v2\n
\pf\, The main part of the proof  will be given by means of  the dual  walk,  denoted by $\hat S$ (defined by $\hat S^x_k=-S_{k}+x$).
If $x\geq \eta\sqrt N$ for some constant $\eta>0$ the assertion is easy to show, since then $P[ S_n^{b_N} =x, n<\zeta ] \asymp P[ S_n^{b_N} =x]\asymp 1/\sqrt N$ (for $x, n$ subject to the condition of the lemma)  according to {\bf P1}(b) and the problem is the same for the bridge without killing.  For the proof of the lemma    it  therefore  suffices to show that if  $\eta < \frac12 \sigma\e$ and  $\tau(\eta) $ is the last  time $\leq n$ when $S^{b_N}$ leaves  $[\eta\sqrt N, \infty)$, then
\beqn\label{N0}
 \sup_{0<x< \eta\sqrt{N}} P\big[ \La_{0,\tau(\eta)+1 }(\de) >  {\textstyle \frac12}\e\, \big|\, S_n^{b_N} =x, n<\zeta \big] \to 0,
 \eeqn
(as $N\to\infty$ and $\de\to 0$) since $ \La_{\tau(\eta) +1, n}(\de) \leq \eta/\sigma < \frac12 \e$.   In order to  separate  the increment $S_{\tau(\eta)+1}-S_{\tau(\eta)}$ from $\La_{0,\tau(\eta) +1}(\de)$ we use the  inequality
$$\La_{0, \tau(\eta) +1}(\de)\leq  \La_{0, \tau(\eta)}(\de) + {\textstyle \frac1{\sigma \sqrt{ N}}}  |S^{b_N}_{\tau(\eta)}-S^{b_N}_{\tau(\eta)+1}|.$$
 Now we switch the description to that by the dual walk.   Write $\tau^x_*(\eta)$ for $\hat \sigma^x_{[\eta\sqrt N,\, \infty)}$ where $\hat \sigma^x$ denotes the hitting time for $\hat S^x$.
 For any  $\e'>0$  we can  choose  $\eta < \frac12 \sigma\e$ and  $N_0>1$ so that  
\beqn\label{eq_l4.2}
\sup_{0<x<\eta \sqrt{N}}P\big[ \,\tau^x_*(\eta) > n/2\;\big|\, \hat S^{x}_n=b_N, n< \hat \sigma^{x}_{(-\infty,0]} \big] <\e' \quad \mbox{for}\;\; N>N_0,
\eeqn  
for, by {\bf P0} and   {\bf P1}(a),  the conditional probability above is dominated by
$$\frac1{p^n_{(-\infty,0]}(x,b_N)}\sum_{z=1}^{\eta\sqrt N}p^{n/2}_{\{0\}}(x,z)p^{n/2}_{\{0\}}(z,b_N) 
  <  C \frac{\eta^2N x}{(n/2)^{3/2}}\cdot \frac{\eta \sqrt N b_N}{(n/2)^{3/2}}\bigg/ \frac{xb_N}{N^{3/2}} < C_{\e}\eta^3.
$$
    With the help of strong Markov property of $\hat S^x$ one  applies   what is mentioned above for the case $x\geq \eta \sqrt{N}$ to the walk $\hat S^x$ to find $\de>0$ so that 
\[
 \sup_{0<x<\eta\sqrt {N}}P\big[ \La_{0,\tau^x_*(\eta)}(\de) > {\textstyle \frac14}\e\, \big|\, S_n^{b_N} =x, n<\zeta \big] < \e'.
 \]
  Since   in  this bound  as well as in (\ref{eq_l4.2}) $\e'$ may be made arbitrarily small, (\ref{N0}) follows  if we can show that
\beqn\label{N1}
 P\bigg[\,  \frac1{\sigma \sqrt{ N}}  \big(\hat S^{x}_{\tau^x_*(\eta)}-\hat S^{x}_{\tau^x_*(\eta)-1}\big) > {\textstyle \frac14}\e, \; \tau^x_*(\eta)\leq {\textstyle \frac12}n\, \bigg |\,   \hat S^{x}_n=b_N, n< \hat \sigma^{x}_{(-\infty,0]}\,\bigg] \to 0
 \eeqn 
 for any prescribed $\eta>0$.

 Our proof of (\ref{N1}) rests on the following facts: for $0<x<\eta\sqrt N$, 
\beqn\label{M1}
\begin{array}{ll}
{\rm (a)}  \qquad  P\big[\, \hat S^{x}_{\tau^x_*(\eta)}- \eta \sqrt{N} >   \eta\sqrt N, \tau^x_*(\eta) <\hat \sigma^x_{(-\infty,0]} \big] =(x/\eta\sqrt{N})\times o(1), \qquad \qquad\\[2mm]
{\rm (b)} \qquad   p^n_{(-\infty,0]}(x, b_N) \geq c_\e\, x/n \quad (\e N<n <N).
\end{array}
\eeqn
(In (b) $c_\e$ is a positive constant that may depend on  $\e$.) Here   (a) follows from   {\bf P5}  
and (b)  from {\bf P0}.
 Let ${\cal E}$ stand for the event  
$${\cal E} := \Big\{ {\textstyle \frac1{\sigma \sqrt{ N}}}  \big(\hat S^{x}_{\tau^x_*(\eta)}-\hat S^{x}_{\tau^x_*(\eta)-1}\big) > {\textstyle \frac14}\e, \tau^x_*(\eta)\leq {\textstyle \frac12}n \wedge \hat \sigma^x_{(-\infty,0]}\Big \}.$$
For simplicity we suppose $n$ to be even.
 Applying  (\ref{M1}b) to the dual walk with the help of the trivial bound  $p^{\lfloor n/2\rfloor}_{(-\infty,0]}(z, b_N) \leq C/\sqrt n$ leads to  
\beqn\label{M2}
P\big[\hat S^{z}_{n/2}=b_N, n/2< \hat \sigma^{z}_{(-\infty,0]}\,\big] \leq C'(\sqrt n/x) P\big[\hat S^{x}_n=b_N, n< \hat \sigma^{x}_{(-\infty,0]}\,\big] \qquad (z>0)
\eeqn
and hence 
 \beq
 P\big[\, {\cal E},   \hat S^{x}_n=b_N, n< \hat \sigma^{x}_{(-\infty,0]}\,\big] 
&=& \sum_{z =1}^\infty P\big[ {\cal E}, \hat S^x_{n/2}=z \big]  P\big[\hat S^{z}_{n/2}=b_N, n/2< \hat \sigma^{z}_{(-\infty,0]}\,\big] \\
 &\leq&    C' (\sqrt n /x)P({\cal E})   P\big[\hat S^{x}_n=b_N, n< \hat \sigma^{x}_{(-\infty,0]}\,\big]. 
\eeq 
The constant $\eta$ may be supposed to be smaller than  $ {\textstyle \frac18}\sigma \e$ so that the occurrence of ${\cal E}$ entails that $\hat S^{x}_{\tau^x_*(\eta)} > \frac14\sigma \e \sqrt N > 2\eta \sqrt{N}$, and hence  that 
$${\cal E}\subset \big \{ \hat S^{x}_{\tau^x_*(\eta)}- \eta \sqrt{N} > \eta \sqrt N, \tau^x_*(\eta)\leq \hat\sigma^x_{(-\infty,0]} \big\},$$
hence  
$P({\cal E}) = [x/\eta\sqrt N]\times o(1)$
 according to (\ref{M1}a), and we may conclude 
$$
 P\big[\, {\cal E},   \hat S^{x}_n=b_N, n< \hat \sigma^{x}_{(-\infty,0]}\,\big] = P\big[\hat S^{x}_n=b_N, n< \hat \sigma^{x}_{(-\infty,0]}\,\big] \times o(1),
 $$
which is the same as (\ref{N1}). Proof of Lemma \ref{lem4.2} is finished. \qed

\v2 \v2
\subsection{Case $E[ |S_1|^3; S_1<0]=\infty$} 

Here we suppose 
in addition to $E[|S_1|^3; S_1]=\infty$ 
that $F(-u)$ is regularly varying as $u\to \infty$ with index $-\beta \in [-3,-2]$.
 
Given  $M > 1$ we shall let $x\wedge (-y)\wedge n\to\infty$ under the constraint
\beqn\label{T5}
   M^{-1}\leq x/\sqrt n \leq M \quad\mbox{and}\quad  M^{-1}\leq -y/\sqrt n \leq M.
   \eeqn
The following result refines the estimates  given in  {\bf P3} and {\bf P4} under the present assumption.
   \begin{Prop}\label{prop2}  As   $x\wedge (-y)\wedge n\to\infty$ under (\ref{T5})
\v2   
   {\rm (i)} \; if $2\leq \beta <3$,
$$ \qquad {\displaystyle    p^n_A(x,y) \sim \sum_{k=1}^n\sum_{z=1}^\infty \sum_{w=1}^\infty \mathfrak{g}_{\sigma^2 k}^0(x,z) \frac{\beta F(-w-z)}{w+z} \mathfrak{g}_{\sigma^2(n-k)}^0(-w,y);}$$

{\rm (ii)} \;  if $ \beta =3$ and  $E\big[|S_1|^3;S_1<0\big] =\infty$,     
$$  p^n_A(x,y) \sim  \rho_{\sigma^2 n}(x+|y|) \frac{2}{\sigma^2}\int_0^{\sqrt n} F(-u)u^2du.$$
\end{Prop}

\v2
We write
$$\xi_N= x/\sqrt{\sigma^2N} \quad \mbox{and}\quad  \eta_N=  y/\sqrt{\sigma^2N}.$$
Recall the definition of $J_t(\xi,\eta)$.  An elementary computation derives from Proposition \ref{prop2}  the following
   \begin{Cor}\label{cor20}   For  any $M > 1$ uniformly  for $x, y, n$ satisfying  (\ref{T5}) and  $M^{-1}<n/N<M$, as  $N\to \infty$,
 $$  p^n_A(x,y) \sim \left\{\begin{array}{ll}  \sigma^{-\beta-1} \big[ \sqrt N F(-\sqrt N\,) \big]J_{n/N}(\xi_N,\eta_N)  \quad  &\beta<3,\\[2mm]
\sigma^{-4}\big [2 \int_0^{\sqrt N} F(-u)u^2du  /  N\big]\rho_{n/N}(\xi_N+\eta_N)  \quad & \beta =3.
 \end{array} \right.
 $$
 \end{Cor}
 
 \v2
The  convergence of finite dimensional distributions of $X^{N}$---given by a formula 
 analogous to  (\ref{T00})---follows from Corollary \ref{cor20}.  The proof of Theorems \ref{thm1} and \ref{thm2} is accordingly  accomplished   if we verify  Proposition \ref{prop2} as well as the tightness of the conditional law
 of $X^{(N)}$. The verification of Proposition \ref{prop2} is given in the paragraphs 4.2.1 and 4.2.2 for the cases  $\beta<3$ and $\beta=3$ respectively.  The tightness proof is given in the paragraph 4.2.3. 

\v2\v2

{\bf 4.2.1.} {\it Proof of {\rm (i)} of Proposition \ref{prop2}.}
\v2

Let $\beta<3$ and $\zeta_x$ and $\zeta_x'$ be defined as in (\ref{T000}) with $x$ in place of $b_N$.  We first observe that 
  as $x\wedge (-y)\wedge n \to\infty$ under (\ref{T5}) and $\e\downarrow 0$ (independently)
 \beqn\label{T6}
 P\big[\,  S^x_{\zeta_x-1} \wedge (- S^x_{\zeta_x}) > \e \sqrt n \,\; \mbox{and} \;\, \zeta'_x= \zeta_x-1\, \big|\, S^x_n =y,\, \sigma_A^x >n\big] \to 1
 \eeqn
 and
 \beqn\label{T60}
 P\big[\, \e n < \zeta_x < (1-\e)n  \, \big|\, S^x_n =y, \sigma_A^x >n\, \big] \to 1.
 \eeqn
 ((\ref{T60}) is proved also in case $\beta =3$  but by a  different approach in the next paragraph {\bf 4.2.2}.)
 (\ref{T6}) follows immediately from  Corollary \ref{cor3.1}(a,b).  
  For the proof of (\ref{T60}), recalling the convention $\max A = 0$ consider the representation 
 \beqn\label{T61}
 p_A^n(x,y)= \sum_{k=1}^n\sum_{z=1}^\infty \sum_{z'\in (-\infty,0]\setminus A} p_{(-\infty,0]}^{k-1}(x,z)  p(z'-z) p_A^{n-k}(z',y).
 \eeqn
 It then suffices to show    that  for any $\eta>0$   there exists $\e >0$ such that the sum
  in (\ref{T61})     restricted to $ k  \leq \e n$ is at most $\eta$, the part $k>(1-\e)n$ being disposed of by duality because of (\ref{T6}). 
 We may restrict the summation over $z$ and $z'$ to  $z\wedge (-z')>\th \sqrt n$ with a positive constant $\th$  because  of (\ref{T6}) again . The proof is therefore finished if we show that for any $\eta>0$ and $\th>0$  there exists  $\e>0$ such that
$$\sum_{ k\leq \e n}\sum_{z>\th \sqrt n}\sum_{w>\th \sqrt n} \frac{P[S^{x}_{k-1}= z]p(-w-z) p^{n-k}_A(-w,y)}{p^n_A(x,y)} \leq \eta.
$$
On using  the trivial bound $p^{n-k}_A(-w,y) \leq C/\sqrt {n-k}$ and  ${\bf P3}$ in turn  this triple sum is dominated by
$$\frac{C}{\sqrt n}\sum_{k\leq \e n} \frac{F(-2\th \sqrt n) 
}{p^n_A(x,y)} \leq C'\th ^{-\beta} \frac{\e n^{1/2} F(-\sqrt n)}{p^n_A(x,y)} \leq \varepsilon C_M\frac{\th^{-\beta}n^{3/2}F(-\sqrt n)}{1+\sum_{w=1}^{\sqrt n} w^3 p(-w)}.
$$
Since  the last ratio is bounded,  the right-most member   becomes arbitrarily small along with  $\e$,
as desired. Thus (\ref{T60}) has been proved.

   In view of  {\bf P0} $p_{(-\infty,0]}^{k}$ and $p_A^{k}$ being interchangeable, 
  we infer from     (\ref{T6}) and (\ref{T60})  that  
 \beq
p_A^n(x,y) = \sum_{k=\lfloor\e n\rfloor}^{\lfloor (1-\e)n\rfloor} \,\sum_{z\geq \e\sqrt n} \,\sum_{w\geq \e \sqrt n} p_A^{k}(x,z)  p(-w-z) p_A^{n-k}(-w,y)\{1+o_\e(1)\},
\eeq
 where $o_\e(1) \to 0$ as  $n\to\infty$ under (\ref{T5}) and $\e\downarrow 0$. In  this   triple sum  we may replace $p^{n-k}_A(-w,y)$ by  $\mathfrak{g}_{\sigma^2(n-k)}^0(-w,y)$ owing to {\bf P1}  since the contribution to the sum 
 from $w>M'\sqrt n$ becomes negligible as $M'\to\infty$ as is assured  without difficulty. 
Similarly   $p^k_A(x,z)$  may be replaced by $\mathfrak{g}_{\sigma^2k}^0(x,z)$. 
 Since either of the sums over $0<w< \e\sqrt n$  and  $0<z< \e\sqrt n$  contributes only $o_\e(1)$ and   for reason of  symmetry  Proposition \ref{prop2} (i)  therefore follows from the next lemma.

\begin{lem}\label{lem4.3} Suppose that $F(-t) = t^{-\beta}L_-(t)$,   $t\geq 1$ with a slowly varying function  $L_-$ and $\beta>1$.
Then for each $M > 1$,  as $z\wedge (-y)\wedge n \to\infty$ under $M^{-1}\leq -y/\sqrt n \leq M $ and $z\leq M\sqrt n$,
 \beqn\label{eq_3} \sum_{w=1}^\infty p(-w-z) \mathfrak{g}_{n}^0(-w,y)\, \sim \, \sum_{w=1}^\infty \frac{\beta L_-(w+z)}{(w+z)^{\beta+1}} \mathfrak{g}_{n}^0(-w,y).
 \eeqn
 \end{lem}
 \v2\n
 \pf\,
 On summing by parts, 
$$
 \sum_{w=1}^\infty p(-w-z) \mathfrak{g}_{n}^0(-w,y) =  \sum_{w=1}^\infty F(-w-z)[ -\mathfrak{g}_{n}^0(-w+1,y) + \mathfrak{g}_{n}^0(-w,y)].
$$
  We can choose a slowly varying function $\tilde L$ so that   $L_-(t) = \tilde L_-(t)(1+\de(t))$ with $\tilde L_-'(t) =o(L_-(t)/t)$ and $\de(t)\to 0$ so that 
  $[t^{-\beta}\tilde L_-(t)]' = (\beta t^{-\beta-1})L_-(t) (-1+o(1))$ as $t\to\infty$. Then
 substituting
$F(-w-z) = (w+z)^{-\beta}\tilde L_-(w+z)(1+\de(w+z))$, and summing by parts back, the sum of the last series above is written as  
\beqn\label{T4}
  \sum_{w=1}^\infty  \beta(w+z)^{-\beta-1}L_-(w+z) \mathfrak{g}_{n}^0(-w,y)\{1+o(1)\} 
+ r(z,y,n),
\eeqn 
 where
$r(z,y,n) =   \sum_{w=1}^\infty [\mathfrak{g}_{n}^0(-w+1,y) - \mathfrak{g}_{n}^0(-w,y)]  (w+z)^{-\beta}\tilde L_-(w+z)\de(w+z)$.
Noting 
$$|\mathfrak{g}_{n}^0(-w+1,y) - \mathfrak{g}_{n}^0(-w,y)|< C/n  \quad \mbox{for all }\quad  w, y,$$  
we  deduce that
$$|r(z,y,n)|\leq \frac{C}{n} \sum_{j\geq z} j^{-\beta}\tilde L_-(j)\de(j) = \frac{L_-(z)}{nz^{\beta-1}}\times o(1).$$ 
On the other hand if   $-y \asymp \sqrt n$ and $1< z=O(\sqrt n)$, then the sum on the right side of (\ref{eq_3}) is bounded from below by a positive multiple of
$$\sum_{z\wedge \sqrt n \leq w\leq 2(z\wedge \sqrt n)}\frac{ L_-(w+z) w|y|}{(w+z)^{\beta+1}n^{3/2}}  \geq C' \frac{L_-(z)}{nz^{\beta-1}}$$
with $C'$  a positive constant, showing $r(z,y,n)$ is negligible. This finishes the proof. \qed

\v2\v2
{\bf 4.2.2.}  {\it Proof of {\rm  (ii)} of Proposition \ref{prop2}.}  
\v2
Let $\beta=3$ so that $F(-u) = u^{-3}L_-(u)$, $u>0$ with a slowly varying $L_-$,  and put
$$L^*_-(u) = \int_1^u \frac{L_-(s)}{s}ds,$$
which is also slowly varying.   Recall  $L^*_-(u)/L_-(u) \to\infty$ as $u\to\infty$,  the fact that differentiates the case $\beta=3$  from the case $\beta<3$ (as is exhibited  by Lemma \ref{lem3.3}).

We follow the proof of  {\bf P2}  given in \cite[pp.702-703]{U1dm} for the special case $A=\{0\}$.
As  therein,  take a small  $\e>0$ and break the sum   in the expression of $p^n_A(x,y)$ given in  (\ref{T61})  into three parts by splitting the range of $k$ according as
$1\leq k <\e n,$ $ \e n\leq k <(1-\e)n$ and $(1-\e)n\leq k\leq n$,
and call the corresponding sum  $I, I\!I$ and $I\!I\!I$.

 On using  the second half of {\bf P6} and  the bound $p^{n-k}_A(w,y) \leq C( |wy|\wedge n)/n^{3/2}$  ($k<\e n$) the part $I$ is dominated by a constant multiple of 
$$\sum_{k\leq \e n} \bigg(\frac{1}{\sqrt k x} \sum_{0<w< \sqrt n}   H^{+\infty}_{(-\infty,0]}(-w)\frac{|wy|}{n^{3/2}}  + \frac1{\sqrt n\,} \sum_{w>\sqrt n} p_A^{k-1}(x,-w) F(-w) \bigg).$$
By the present assumption  about  $F$
\beqn\label{H/ww}
H^{+\infty}_{(-\infty,0]}(-w) \sim \frac{2}{\sigma^2} \sum_{z>w}F(-z) \sim \frac{L_-(-w)}{\sigma^2 w^{2}}\quad (w\to \infty)
\eeqn
(cf. \cite[(2.7)]{U1dm.f.s}). Under  (\ref{T5}) $y/x \leq M^2$ and    easy computations  deduce  that 
\beqn\label{I}
I \leq C_M[\sqrt{\e}\, L^*_-(\sqrt n\, )/n +\e L_-(\sqrt n\,)/n] <C_M'\sqrt\e \, L^*_-(\sqrt n\,)/n. 
\eeqn
$I\!I\!I$ admits the same bound  as a dual relation since $P^*_{(c,N)}[\zeta=\zeta'+1] \to 1$ by Corollary \ref{cor3.1}.

As for $I\!I$ we note that for arbitrarily small $\eta>0$  the range of the variable $w$ may be restricted to  $w<\eta\sqrt n$ in view of Lemma \ref{lem3.3}(c).
Then 
by the first  half of {\bf P6}  and {\bf P1}(d)
$$I\!I=\sum_{ \e n\leq k \le (1-\e)n}\frac{f_+(x)\mathfrak{g}_{\sigma^2 k}(x)}{ k}\,\sum_{0<w<\eta \sqrt n}H^{+\infty}_{(-\infty,0]}(-w) p_A^{n-k}(-w,y)(1+o_{\e,\eta}(1)).$$
Here  for each pair of $\e$ and $\eta$,  $o_{\e,\eta}(1)\to 0$ as $n\to\infty$  uniformly for $x,y$ satisfying (\ref{T5}).

By {\bf P1} $p_A^{n-k}(-w,y) \sim 2[\sigma^{2}(n-k)]^{-1} |wy|\mathfrak{g}_{\sigma^2(n-k)}(y)$ as $w\to\infty$ under $w=o(\sqrt n)$  and observe,  on  replacing $f_+(x)$   by $x$, 
\[
I\!I=\sum_{ \e n\leq k \le (1-\e)n}\frac{ 2x|y|\mathfrak{g}_{\sigma^2k}(x)\mathfrak{g}_{\sigma^2(n-k)}(y)}{\sigma^2 k(n-k)}\,\sum_{0<w<\eta \sqrt n} H^{+\infty}_{(-\infty,0]} (-w)w(1+o_{\e,\eta}(1)). \nonumber
\]
 Writing  $x_n = x/\sigma \sqrt{n}$,  $y_n = y/\sigma \sqrt{n}$ we see
 $x\mathfrak{g}_{\sigma^2k}(x)/  k= \rho_{k}(x/\sigma) = \rho_{k/n}(x_n)/n$  and 
\beq
\sum_{ \e n\leq k \le (1-\e)n}\frac{2x|y| \mathfrak{g}_{\sigma^2k}(x)\mathfrak{g}_{\sigma^2(n-k)}(y)}{\sigma^2 k(n-k)}
= \frac2{n\sigma^2}\int_0^1 \rho_t(x_n) \rho_{1-t}(|y_n|)dt + O\bigg(\frac{\e}{n}\bigg) +o\bigg( \frac{1}{n}\bigg).
\eeq
Here we have used the condition (\ref{T5}) again. On the other hand  in view of   (\ref{H/ww})
$$\sum_{0<w<\eta \sqrt n} H^{+\infty}_{(-\infty,0]} (-w)w(1+o_{\e,\eta}(1)) = \frac1{\sigma^2}L^*_-(\sqrt n\,)(1+o_\e(1)).$$ Finally observing 
$\int_0^1 \rho_t(x_n) \rho_{1-t}(|y_n|)dt =\rho_1(x_n+|y_n|)= \rho_{\sigma^2 n}(x+|y|)\sigma^2n$
we find that
\beqn\label{II}
I\!I= \frac2{\sigma^2} \Big[( \rho_{\sigma^2 n}(x+|y|) + O(\e/n)\Big] L^*_-(\sqrt n\,)(1+o_\e(1)),
\eeqn
which together with the bounds of $I$ and $I\!I\!I$ verifies the formula   in  (ii).
\qed

\v2\v2
   
{\bf 4.2.3.} {\it Proof of  tightness}.   
\v2
  It suffices to show that for any $\e$ there exists $N_0$ and $\de>0$ such that
$$P^*_{(c,N)}[ \La_{0, \zeta-1}(\de) >\e] <\e   \quad \mbox{for}\quad N\geq N_0$$
($\zeta$ and  $\zeta'$ are given in (\ref{T000})), since  $P^*_{(c,N)}[\zeta=\zeta'+1]\to 1$ in view of Corollary \ref{cor3.1}(a).  Observe that $P[\La_{0,\zeta-1}(\de)>\e,   \zeta =k, S^{b_N}_{k} = -w]$, $w\geq 1$ may be written as
$$\sum_{x=1}^\infty P[\La_{0, k-1}(\de) >\e, k\leq \zeta, S^{b_N}_{k-1}= x]p(-w-x).$$
By the dual assertion of  Lemma \ref{lem3.3}(b,c)  one can choose $M$ so large that 
  $P^*_{(c,N)}[ S^{b_N}_{\zeta-1} \geq M\sqrt N\,]<\e/2$ and accordingly  obtains that   $P^*_{(c,N)}[\La_{0,\zeta-1}(\de) >\e]$   is at most $\e/2$ plus
\beqn\label{T7}
\sum_{k=1}^N\sum_{0<x<M \sqrt N}\sum_{w=1}^\infty \frac{P[\La_{0,k-1}(\de) >\e,k\leq\zeta,  S^{b_N}_{k-1}= x]p(-w-x) p^{N-k}_{A}(-w,-c_N)}{p^N_A(b_N,-c_N)}.
\eeqn
By Lemma \ref{lem4.2},   as $\de\downarrow 0$ and $N\to\infty$,   $P[\La_{0,k-1}(\de)> \e\,|\, k\leq \zeta, S^{b_N}_{k-1} =x]\to 0$ or, what amounts to the same, 
\[
P[\La_{0,k-1}(\de) >\e,  k\leq \zeta, S^{b_N}_{k-1} = x] = p^{k-1}_{(-\infty,0]}(b_N,x)\times o(1)
 \]
uniformly for $\eta N \leq k\leq (1-\eta) N$ and $0<x <M\sqrt N$ for each positive $\eta$. 
On recalling (\ref{T61}) this disposes of the contribution from $\eta  N \leq k\leq (1-\eta)N$.
 It follows from (\ref{T60}) (if $\beta<3$)  and (\ref{I}) (if $\beta=3$)  that  the contribution  from   $k<\eta N$  
 is negligible; and similarly for that from $k> (1-\eta)N$, the proof of tightness is complete. \qed

 \v2 
 \subsection{Notes on the proofs of Theorems \ref{thm3} and \ref{thm4}}
 \v2
 
 The proofs are similar to those of Theorems \ref{thm1} and \ref{thm2} and we do not present them but indicate some points that make difference  from the latter. 
Throughout this section we assume  (\ref{stbl0}) to be valid, put $b_N= \lfloor b\la_n\rfloor$ and $c_N= \lfloor c\la_n\rfloor$  and let $P^*_{(c,N)}$, $\zeta$ and $\zeta'$ be defined as before.
 First we note that
 the propositions corresponding to {\bf P1} through {\bf P4} and {\bf P6} for the  case $\sigma^2=\infty$ are obtained
in \cite{Uattrc}: specifically the corresponding results of  \cite{Uattrc} are  
\v2
{\bf P1}(a,b) $\mapsto$ Theorem 3, Corollary 3, Proposition 2.1; \quad  {\bf P1}(c) $\mapsto$ Lemma 5.2 [(\ref{p_p/y})];

  {\bf P1}(d) $\mapsto$ Proposition 2.3(ii)  [(\ref{eqP2})];  \quad  {\bf P2} $\mapsto$  Theorem 4;

{\bf P3} $\mapsto$  Proposition 6.1; \quad  {\bf P4} $\mapsto$  Proposition 2.2; \quad  {\bf P6}  $\mapsto$ Lemma 6.5,
\v2\n
where Theorems,  Propositions or  Lemmas are those of \cite{Uattrc}, of which a few   are  already used in Section 4, their reference numbers being  indicated in the square brackets. 
The result corresponding  to {\bf P0} (in particular  its first one saying $p^n_{\{0\}}(x,y)\sim p^n_{(-\infty,0]}(x,y)$) that is  missing in the above list  and  used in the proofs of Corollary \ref{cor3.1} and Lemma \ref{lem4.2} does not hold in the relevant range  $0< x, y \leq \lambda_n$ but is valid if restricted to  $x,y \in [M^{-1} \lambda_n, M \lambda_n]$.  Corollary \ref{cor3.1}  was used only under this restriction, while  a result corresponding to   Lemma \ref{lem4.2} will be proved below.   
As for  {\bf P5} also missing in the above,  we need a corresponding one in a dual form, which  will be 
presented below  in  (\ref{eqP5}). 

 It is shown in \cite[Theorem 6]{Uattrc} that as $N\to\infty$
\beqn\label{eqU1}
p_A^N(b^N,-c^N) \sim C_A^+ (b+c) \mathfrak{p}_1(b+c)/N \quad \mbox{if}\quad C^+<\infty \eeqn 
where $C_A^+$ is some positive constant.  Instead of Lemmas \ref{lem4.1}  and \ref{lem4.2} we obtain the following Lemmas \ref{lem4.4} and \ref{lem4.5}, respectively.
 \begin{lem}\label{lem4.4} \;  Suppose $C^+<\infty$. Then  there exists a constant $C$ depending on $b, c, A$ and $F$ such that for $R>1$
  \beqn\label{T20}
P^*_{(c,N)}[ \, S^{b_N}_k>  R \; \mbox{for some} \;  \zeta\leq k\leq N \,] \leq C/a(-R)
\eeqn
and $P^*_{(c,N)}[ \,  S^{b_N}_k< -R \; \mbox{for some} \; k < \zeta' \,] \leq C/a(-R)$.
 \end{lem}
\v2\n
\pf\,  As in the proof of Lemma \ref{lem4.1} we see that   $P[\sigma^z_{[R,\infty)} <\sigma_{\{0\}}^z]\leq C_2a(z)/a(-R)$ ($z<0$)  and owing to Lemma \ref{lem3.4} 
$\sup_{x>R} p^{N-k}_{\{0\}}(x,-c_N) \leq C_1 L(c_N)/[c_N]^\a \sim C_1/c^\a N$ ($k<N$) and that  
by (\ref{eqU1}) and (\ref{eqL3.3}) the probability in (\ref{T20}) is at most a constant multiple of 
$$\frac1{p_A^N(b^N,-c^N)} \sum_{z=-\infty}^{-1} P[S^{b_N}_{\zeta} =z] \frac{a(z)}{Na(-R)}\sim C_{A,b,c}  \frac{a(b_N)}{a(-R)} \leq \frac{C}{a(-R)},$$
where we have also used  relation  (\ref{H/aa}) that gives  $E[a(S^{b_N}_{\zeta})]=a(b_N)$. \qed 
\v2

For $0\leq j<k\leq NT$ let $\La''_{j,k}(\de)$ denote  the random variable 
$$\La''_{j,k}(\de) = \sup_{\frac1{N}j \leq s<t\leq \frac1{N}k ;\, |t-s|<\de}\,  \sup_{s<u<t} |X^{(N)}_u-X^{(N)}_s|\wedge  |X^{(N)}_t-X^{(N)}_u|.$$

\begin{lem}\label{lem4.5} For each $0<\e <1/2$, as $N\to\infty$ and $\de\downarrow 0$ independently
  $$ P[\, \La''_{0,n}(\de) > \e\, |\, S_n^{b_N} =x, n<\zeta \,] \to 0 $$
 uniformly for  $\e N <n< N$ and   $0< x < \e^{-1} \la_N$.
\end{lem}
\v2\n
\pf\, As in the proof of Lemma \ref{lem4.2} the proof is  reduced to verification of   (\ref{N1}) modified in an obvious way (except for the existence of $\eta$ that makes   (\ref{eq_l4.2}) valid (with $\eta\sqrt N$ replaced by $\eta\la_N$) for which we use the convergence of the normalized walk conditioned to stay positive to a stable meander (as found in \cite{D1}). 
The verification of  (\ref{N1})  is carried out in the same line  if an appropriate substitute for (\ref{M1}) (and accordingly  that for (\ref{M2})) is
given.  
Let $U_{{\rm as}}$ denote the renewal function of the ascending ladder height process of the walk $(S_n)$. Then,     as $R\to\infty$ 
\beqn\label{eqP5}
P[\hat S^x_{\sigma_{[R,\infty)} }\geq  R+r, \hat\sigma_{[R,\infty)} < \hat\sigma_{(-\infty,0]}] =  \frac{U_{{\rm as}}(x)}{U_{{\rm as}}(R+r) - U_{{\rm as}}(R)}\times o(1)
\eeqn
uniformly for $0< x <R$, $r\geq 1$ according to  \cite[Lemma 1.1]{Uladd}.  From Proposition 11 of \cite{D} we infer   that for each  $M>1$ there exists a constant $c_\circ >0$ such that  
\beq
 p^n_{(-\infty,0]}(y,x) &=& P[\, \hat S^{x}_{n}=y, n <\hat \sigma^x_{(-\infty,0]} ] \\
 &\geq& c_\circ \, U_{{\rm as}}(x)P[\sigma_{[1,\infty)}^0 >n ]/\la_n  \qquad (0< x< M^{-1}\la_n <y < M\la_n).
\eeq
 Moreover $U_{{\rm as}}$ varies regularly (with index $\alpha-1$)  and $U_{{\rm as}}(\la_n)P[\sigma_{[1,\infty)}^0 >n ]$  approaches a positive constant as $n\to\infty$. From these facts one deduces that 
for $0<x<\eta\la_N$, 
\[
\begin{array}{ll}
{\rm (a)}  \qquad  P[\, \hat S^{x}_{\tau_*(\eta)}- \eta \la_N >   \eta \la_N, \tau^x_*(\eta) <\hat \sigma^x_{(-\infty,0]} ] =[U_{{\rm as}}(x)/U_{{\rm as}}(\la_N)]\times o(1), \qquad \qquad\\[2mm]
{\rm (b)} \qquad   P[\, \hat S^{x}_{n}=b_N, n <\hat \sigma^x_{(-\infty,0]} ]  \geq c_\e\, [U_{{\rm as}}(x)/U_{{\rm as}}(\la_N)]/\la_N \quad (\e N<n <N)
\end{array}
\]
instead of (\ref{M1}), and then 
\[
P[\hat S^{z}_{n/2}=b_N, n/2< \hat \sigma^{z}_{(-\infty,0]}\,] \leq C[U_{{\rm as}}(\la_N)/U_{{\rm as}}(x)] P[\hat S^{x}_n=b_N, n< \hat \sigma^{x}_{(-\infty,0]}\,] \qquad (z>0)
\]
instead of (\ref{M2}). The rest of proof  is easy and omitted. \qed
\v2

With these two lemmas as well as those corresponding to {\bf P1} to {\bf P6} one can follow the arguments of Sections 4.1 and 4.2 to prove Theorems \ref{thm3} and \ref{thm4}.  

\section{Appendix}

(A) \, Let $\ell (t)$ be the local time  at the origin of the Brownian motion $W$.  For a constant $\la>0$ let $\xi^{(\la)}_t,\,t\geq 0$ be
$W$ killed at  the rate $\la \ell (t)$; in other words,  $\xi^{(\la)}$ is a Markov process whose sample path is  continuous and whose  transition probability is given by
$$q^{(\la)}_t(\xi,\eta) = E_\xi[\, e^{-\la \ell (t)}; W_t\in d\eta]/d\eta.$$
Here $P_\xi$ denotes the probability law of  $W^\xi$ and $\xi$ being  dropped from $W^\xi$ under $P_\xi$. Then 
 {\it the conditional law of the process 
$\xi^{(\la)}$ given  $\xi^{(\la)}_0 =b$ and $\xi^{(\la)}_T=-c$ converges weakly, as $\la\to\infty$,  to the law of the  process $X^{b,c,T}$ described in (\ref{X0}) and the  density of  finite dimensional distribution of the limit process given 
on the right side of 
(\ref{T00}) is  expressed as}
\beqn\label{Lim}
\lim_{\la\to\infty} \frac{\la E_b[e^{-\la \ell(T)};  \;\cap_{j=1}^m 
\{W_{t_j}\in d\xi_j \} \; \mbox{and} \;
 \cap_{k=0}^{m'} \{W_{t_j}\in d\eta_j\} ]}{\rho_T(b+c) \, d\xi_1\cdots d\xi_m\cdot d\eta_0 \cdots d\eta_{m'}}
  \eeqn

We prove only  the latter half of the assertion, the other half being argued  similarly to what is done in Section 4.1 for the random walk. The proof rests on   the formula
\beqn\label{IM}
P_0[ \ell(t)\in du,  W_t\in d\eta] =  \rho_t(|\eta| +u)dud\eta
\eeqn
(cf.  \cite[Problem 2.2.3]{IM}, \cite[Eq(VI.2.18)]{RY} etc.), from which  we derive 
$$
\begin{array}{ll}
{\rm (i)} \quad  \lim_{\la \to\infty} q^{(\la)}_t(\xi,\eta) = \mathfrak{g}_t^0(\xi,\eta) \quad &\mbox{if} \quad \xi\eta >0\\[2mm]
{\rm (ii)} \quad \; \lim_{\la \to\infty}\la q^{(\la)}_t(\xi,\eta) = \rho_t(\xi+|\eta|) \quad& \mbox{if} \quad \eta <0<\xi, \qquad\qquad\qquad\qquad\qquad\qquad
\end{array}
$$
as follows. If $\xi \eta>0$, then as $\la\to \infty$, $E_\xi[e^{-\la \ell(t)}; W_t\in d\eta, \sigma^W_{\{0\}} <t]/d\eta \to 0$, so that 
 $$q^{(\la)}_t(\xi,\eta) = E_\xi[e^{-\la \ell(t)}; W_t\in d\eta, \sigma^W_{\{0\}} >t] /d\eta +o(1)
= \mathfrak{g}_t^0(\xi,\eta)\{1+o(1)\},$$
showing (i).
If $\eta <0<\xi$,
(\ref{IM}) together with $\int_0^t \rho_{s}(\xi)\rho_{t-s}(|\eta|+u)ds = \rho_t(\xi + |\eta|+u)$  leads to
\beq 
q^{(\la)}_t(\xi,\eta) &=& \int_0^t \rho_{s}(\xi)E_0[e^{-\la \ell(t-s)}; \; W_{t-s}\in d\eta]ds/d\eta \\
&=&\int_0^\infty \rho_t(\xi+|\eta|+u) e^{-\la u} du,
\eeq
hence (ii) follows.

Now for $\eta, \zeta< 0<\xi$ and $0<s<t$,
\beq
\la E_\xi\big[e^{-\la \ell(t)}; \, W_s \in d\eta,  W_t \in d \zeta \big]  
&=&\la E_\xi \big[e^{-\la \ell(s)} \1_{d\eta}(W_s) E_{W_s}[ e^{-\la \ell(t-s)}; \, W_{t-s} \in d\zeta] \,\big]\\
&\sim& \rho_s(\xi+|\eta|)\mathfrak{g}^0_{t-s}(\eta,\zeta)d\eta d\zeta,
\eeq
disposing of the case  $0=m <m'=1$ (by taking $\xi=b$, $\zeta=-c$,  and  on using this result   for  $-c= \zeta, \eta< 0<\xi $ and $0<s<t <T$,
\beq
&&\la E_b\big[e^{-\la \ell(T)}; \, W_s \in d\xi, W_t\in d\eta,  W_T\in d\zeta \big]/d\xi d\eta d\zeta  \\
&&=\la E_b\big[e^{-\la \ell(s)} \1_{d\xi}(W_s) 
E_{W_s} [e^{-\la \ell(T-t)}; \, W_{t-s} \in d\eta,  W_{T-t} \in d \zeta] \,\big]/d\xi d\eta d\zeta  \\
&&\longrightarrow\;  \mathfrak{g}_s^0(b,\xi) \rho_{t-s}(\xi+|\eta|)\mathfrak{g}^0_{T-t}(\eta,-c)
 \eeq
($\la\to\infty$), showing  the limit in  (\ref{Lim}) to  agree with  the corresponding density in (\ref{T00})  if  $m=m'=1$. The  case $m'<m=1$ is dealt with in a similar way and   the  general case    by  the double induction on  $(m, m')$.

 \v2\n
 (B) \, Here is given a proof of (\ref{ent_l}). The proof is the same as for the Brownian case but we need to take care of  asymmetry of the processes $Y^\xi$. Let $\xi>0$. 
Let $\sigma_0^{\xi,Y}$ be the first hitting time of zero for $Y^\xi$.
 Noting that  the event  $\sigma^{\xi,Y}_0>s$  entails $Y^\xi_s >0$ since $Y^\xi$ makes no negative jumps and $\mathfrak{p}^{\{0\}}_s(\xi, \eta) =\mathfrak{p}^{\{0\}}_s(-\eta,-\xi)$, we see
 \begin{eqnarray*}
 P[\sigma^{\xi,Y}_0>s+t] &=&\int_0^\infty P[\sigma_0^{\xi,Y}>s, Y^\xi_s\in d\eta] P[\sigma^{\eta,Y}_0>t]\\
&=&\int_0^\infty \mathfrak{p}^{\{0\}}_s(\xi,\eta) P[\sigma_0^{\eta,Y}>t]d\eta\\
&=&\int_t^\infty du  \int_0^{\infty}\rho_u(\eta)\mathfrak{p}^{\{0\}}_s(-\eta,-\xi)d\eta,
\end{eqnarray*}
and  differentiation yields $\rho_{s+t}(- \xi)= \int_0^ {-\infty}  
\rho_t(\eta) \mathfrak{p}_s^{ \{0 \}}(\eta,- \xi)d \eta$ that is the same as (\ref{ent_l}).

\end{document}